\PassOptionsToPackage{numbers}{natbib}
 \documentclass[onecolumn,nopreprint]{jasatex}
\pdfoutput=1
\usepackage[breaklinks,colorlinks=true,linkcolor=blue,citecolor=red,backref=page]{hyperref}

\usepackage[utf8]{inputenc}
\usepackage[T1]{fontenc}

\usepackage{amsmath}
\usepackage{amssymb}
\usepackage{amsfonts}
\usepackage{graphicx}
\usepackage{color}

\newcommand{\eps}{\varepsilon} 
 
\newcommand{\RR}{{\mathbb R}} 
\newcommand{\EE}{{\mathbb E}}

\newcommand{\CC}{{\mathbb C}}

\DeclareMathAlphabet{\itbf}{OML}{cmm}{b}{it}

\definecolor{darkgreen}{rgb}{0.1,0.7,0.1}
\definecolor{darkred}{rgb}{0.7,0.1,0.1}

\begin{document}

\title{Acoustic and geoacoustic inverse problems in randomly perturbed shallow-water environments}

\author{Laure Dumaz}
\affiliation{Ceremade, Universit\'e Paris-Dauphine,
75016 Paris - France,
and\\
Sivienn, 
29570 Roscanvel - France}

\author{Josselin Garnier}
\affiliation{Centre de Math\'ematiques Appliqu\'ees, Ecole Polytechnique, 91128 Palaiseau Cedex - France, 
and\\
Sivienn, 29570 Roscanvel - France}

\author{Guilhem Lepoultier}
\affiliation{Sivienn, 29570 Roscanvel - France}

\date{\today}

\begin{abstract}
The main goal of this paper is to estimate the regional acoustic and geoacoustic shallow-water environment 
from data collected by a vertical hydrophone array and transmitted by distant time-harmonic point sources.
We aim at estimating the statistical properties of the random fluctuations of the index of refraction
in the water column and the characteristics of the sea bottom.
We first explain from first principles how acoustic wave propagation can be expressed as Markovian dynamics for 
the complex mode amplitudes of the sound pressure, which makes it possible to express the cross moments of the sound pressure
in terms of the parameters to be estimated.
We then show how the estimation problem can be formulated as a nonlinear inverse problem 
using this formulation, that can be solved by minimization of a misfit function. 
We apply this method to experimental data collected by the ALMA system (Acoustic Laboratory for 
Marine Applications).
\end{abstract}



\maketitle

\section{Introduction}

In this paper we consider acoustic wave propagation in a randomly perturbed
shallow-water waveguide with an absorbing sea bottom.
The random perturbations of the index of refraction are due to internal waves,
which induce  temperature  and  salinity  fluctuations, 
and the sea bottom  is made of sediments, which induce dissipation.
In the regime where random perturbations and dissipation are small and propagation distance is large,
it is possible to get an effective description of the acoustic wave propagation
in terms of a Markovian dynamics for the complex mode amplitudes of the expansion of the pressure field 
onto the guided modes of the unperturbed and non-dissipative waveguide.
This Markovian dynamics involves coupling terms between 
guided modes and mode-dependent dispersion and loss terms.
The coupling terms between guided modes come from the random perturbations in the water column.
The effective dispersion and loss terms come from two effects: the random perturbations in the water  column
induce coupling between the guided and radiative modes 
and the deterministic dissipation in the sediment layer induce an exponential decay of the guided mode amplitudes.
Both effects generate an irreversible, mode-dependent loss of energy carried by the guided modes.

The mathematical literature contains a lot of results on wave propagation in randomly perturbed waveguides 
motivated by underwater acoustics \cite{kohler77,GP07,gomez}.
Those results derive from first principles and make it possible to relate the coefficients of the effective Markovian model
to the physical parameters of the waveguide (in particular, the statistics of the fluctuations of the index 
of refraction and the complex acoustic impedance of the sea bottom).
The mathematical statements also make it clear in which sense the Markovian model approximates the random dynamics
of the mode amplitudes in the waveguide,
because there are subtle effects that follow from the fact that the results are first established in a weak topology, 
as explained in Appendix \ref{app:A}.
Comparisons between theory and numerical simulations show good agreement, 
both for direct and synthetic inverse problems \cite{borcea15},
but there is not so far much comparison between the detailed mathematical predictions and real experiments.
On the other hand, the physical literature contains many theoretical results essentially based on coupled mode equations \cite{beran,dozier,dozier2,creamer,colosi09,colosi12}.
However, the equations are derived ad hoc and it is not always straightforward to relate the coefficients
of the effective equations to the physical
and statistical parameters of the medium. 
This is problematic as we have in mind to use such results to solve an inverse problem.

As we will see in this paper, the effective attenuation of the mode amplitudes plays a key role.
In the physical literature, modal attenuation coefficients are introduced directly in the coupled mode equations without derivation from first principles
 \cite{dozier83,creamer,jensen}.
On the other hand most of mathematical studies do not take into account the attenuation of the bottom layer.
However, even though dissipation is weak, it plays an important role in long-range propagation.
For the direct problem point of view, the equipartition regime that is well characterized by Gaussian statistics \cite{flatte}
and a scintillation index (relative intensity variance) 
close to one in the absence of dissipation looks very different in the presence of weak attenuation and it
may give rise to high fluctuations of the intensity. This was first pointed out by Creamer \cite{creamer}.
For the inverse problem point of view, signals measured on an hydrophone array can be processed
to extract these attenuation coefficients, which in turn allow to get information on the medium,
in particular, the bottom properties, as we will see in Section~\ref{sec:ip}.
We exploit data from at-sea experiments based on DGA's ALMA (Acoustic Laboratory for 
Marine Applications) system \cite{alma15}.
We use data collected by a system consisting of four 32-element, vertical-line hydrophone array and a moored pinger located far away (at 9 km)
and  transmitting time-harmonic signals between $2$ and $13~kHz$ with a three-minute repetition rate.
The measurements were carried out in November (7th-17th) 2016 near the shores of North East Corsica \cite{alma17c}.
The signals recorded by the hydrophones can be processed by cross correlation calculations
to estimate the properties of the medium. We explain this procedure  and report its results  in this paper.

The paper is organized as follows.
The direct problem is analyzed in Sections \ref{sec:intro}-\ref{sec:randomwaveguide}. Section \ref{sec:homo} is a review of the modal decomposition
of the sound pressure in a homogeneous, non-dissipative waveguide. Section \ref{sec:randomwaveguide} describes the Markovian dynamics
of the mode amplitudes in a random, dissipative waveguide.
The inverse problem is formulated and solved using experimental data in Section \ref{sec:ip}.

\section{Wave propagation in waveguides}
\label{sec:intro}%
Our model consists of a two-dimensional waveguide with range axis denoted by $x \in \RR$ and transverse coordinate
denoted by $z \in [0,+\infty)$. We suppose that the depth of the sea is constant and equals $z_{\rm b}$. When $z \in [0,z_{\rm b}]$ the medium is water, when $z> z_{\rm b}$ it becomes sediments. A point-like source at a fixed position $(x,z) = (0,z_0)$ transmits a time-harmonic signal at frequency $\omega$ 
which is collected by a vertical array of receivers (hydrophones) at $x=x_{\rm a}$.

The Helmholtz equation for the acoustic pressure $\hat{p}(x,z)$ writes:
\begin{align}
\Big[  
 \rho(x,z) \nabla \cdot \Big( \frac{1}{\rho(x,z)} \nabla\Big) + \frac{\omega^2}{c(x,z)^2}  \Big] 
\hat{p}(x,z) =  \delta(x)\delta(z-z_0),   \label{eq:pressure0}
\end{align}
for $x \in \RR, \; z \geq 0,$
where $\nabla=(\partial_x,\partial_z)$, $c(x,z)$, resp. $\rho(x,z)$, is the speed, resp. the density, at position $(x,z)$.  
We consider that $\rho$ is stepwise constant and equal to $\rho_{\rm w}$  in the water and 
$\rho_{\rm s}$ in the sediments.
Therefore we have
\begin{align}
\Big[  (\partial_x^2 + \partial_z^2 ) +  \frac{\omega^2}{c(x,z)^2}   \Big] 
\hat{p}(x,z)  = \delta(x)\delta(z-z_0),  \label{eq:pressure}
\end{align}
for $ x \in \RR, \; z \in (0,z_{\rm b})\cup(z_{\rm b},+\infty)$.

{\bf Remark.}
{\it
The model (\ref{eq:pressure}) can be derived from a more realistic three-dimensional situation, in which 
the sound pressure field $\hat{P}$ satisfies in cylindrical coordinates 
$$
\Big[ \partial_r^2  +\frac{1}{r} \partial_r +\frac{1}{r^2} \partial_\theta^2  +\partial_z^2  +
\frac{\omega^2}{c(r,z)^2} \Big] 
\hat{P} = \frac{1}{2\pi r} \delta(r)\delta(z-z_0) .
$$
The solution is radially symmetric and, 
if we neglect a near field factor of the form $\hat{p}/r^{5/2}$ then the scaled pressure field 
$\hat{p}(r,z) = \sqrt{r} \hat{P}(r,z)$
satisfies (\ref{eq:pressure}).}

The differential equation \eqref{eq:pressure} is endowed with the following boundary conditions: 
\begin{equation}
\begin{array}{l}
\bullet \mbox{ Dirichlet boundary condition on the top of the water column: }
\mbox{$\quad \hat{p}(x,0) = 0$ for all  $x \in \RR$,}\\
\bullet \mbox{  Continuity at depth $z_{\rm b}$:}
\mbox{$\quad \hat{p}(x,z_{\rm b}^-) = \hat{p}(x,z_{\rm b}^+)$  and $\partial_z \hat{p}(x,z_{\rm b}^-)/\rho_{\rm w} =   \partial_z \hat{p}(x,z_{\rm b}^+)/\rho_{\rm s}$.}
\end{array}
\label{eq:bc0}
\end{equation}

We are interested in solutions of (\ref{eq:pressure}) such that
$$
\hat{p} {\bf 1}_{(0,+\infty)}(x) 
\in {\cal C}^0\big( (0,+\infty),H_0^1(0,+\infty)\cap H^2(0,+\infty)\big)
\cap  {\cal C}^2\big( (0,+\infty),L^2(0,+\infty)\big) .
$$

Let us now describe the results we obtain in this model. 
We compute in Section \ref{subsec:idealwaveguide} the pressure 
for an ideal (homogeneous) waveguide in an exact form. 
We analyze how it behaves when the sound speed in the water is no longer constant but perturbed 
by some small noise and the sediments are weakly dissipative. 
The wave modes interact with each other and an asymptotic regime in the small noise and large distance limit is studied in Section \ref{sec:randomwaveguide}. 
In particular we compute correlations of the recorded pressure signals for hydrophones 
located in a vertical segment  for exponentially decaying correlations of the medium in Section \ref{subsec:verticalcorrelations}.

\section{Homogeneous waveguide}\label{subsec:idealwaveguide}
\label{sec:homo}%
In this section, we consider a wave speed $c(x,z)$ which is constant both in the water and in the sediments:
\begin{align}
c_0(z) = c_{\rm w} {\bf 1}_{[0,z_{\rm b}]}(z) + c_{\rm s} {\bf 1}_{( z_{\rm b},+\infty)} (z) ,
\end{align} 
with $c_{\rm s} > c_{\rm w}$ (Pekeris waveguide). There is no dissipation, no fluctuation along the $x$-axis.
The analysis of the perfect waveguide is classical \cite{wilcox}
but we include it here for the sake of completeness.
Let us  introduce the Helmholtz operator
\begin{align}
{\cal H} = 
\rho_0(z) \partial_z \rho_0(z)^{-1} \partial_z + \frac{\omega^2}{c_0(z)^2}  ,
 \label{eq:helmholtz}
\end{align}
where $\rho_0(z) = \rho_{\rm w} {\bf 1}_{[0,z_{\rm b}]}(z) + \rho_{\rm s}  {\bf 1}_{( z_{\rm b},+\infty)} (z)$,
with boundary conditions:
\begin{equation}
\label{eq:bc}
\begin{array}{l}
\bullet \mbox{ Dirichlet at the top: $\phi(0) =0$,}\\
\bullet \mbox{  Continuity at $z = z_{\rm b}$: $\phi(z_{\rm b}^+) = \phi(z_{\rm b}^-)$ and 
$\phi'(z_{\rm b}^-)/\rho_{\rm w} = \phi'(z_{\rm b}^+)/\rho_{\rm s}$.}
\end{array}
\end{equation}
We also denote $k_0(z) = k_{\rm w} {\bf 1}_{[0,z_{\rm b}]}(z) + k_{\rm s} {\bf 1}_{( z_{\rm b},+\infty)}(z)$,
with $k_{\rm w} = \omega/c_{\rm w}$ and $k_{\rm s} = \omega/c_{\rm s}$.
The Helmholtz operator  is self-adjoint with respect to the scalar product defined on $L^2(\RR^+,\rho_0^{-1})$ by:
\begin{align}
\nonumber
(\phi_1,\phi_2) &:= \int_0^{\infty}  \rho_0^{-1}(z)\; \overline{\phi_1(z)} {\phi_2(z)}dz \\
&= \rho_{\rm w}^{-1}  
\int_0^{z_{\rm b}} \overline{\phi_1(z)} {\phi_2(z)}dz +  \rho_{\rm s}^{-1} 
\int_{z_{\rm b}}^{+\infty} \overline{\phi_1(z)} {\phi_2(z)}dz.
\label{scalarproduct}
\end{align}
The Helmholtz operator has a spectrum of the form
\begin{equation}
(-\infty,  k_{\rm s}^2) \cup \{ \beta_N^2,\ldots,\beta_1^2\}  ,
\end{equation}
where the $N$ modal wavenumbers $\beta_j$ are positive and 
$k_{\rm s}^2<\beta^2_N < \cdots <\beta_1^2 <k_{\rm w}^2$.

{\bf Discrete spectrum.}
The $j$th eigenvector associated to the eigenvalue $\beta_j^2$ is
\begin{align}
\phi_j(z) = 
\left\{
\begin{array}{l}
A_j \sin (\sigma_j z/z_{\rm b}) \mbox{ if } 0\leq z\leq z_{\rm b}\\
A_j \sin (\sigma_j ) \exp(-\zeta_j (z-z_{\rm b})/z_{\rm b})\mbox{ if }  z\geq z_{\rm b}
\end{array}
\right.
\end{align}
where 
\begin{align}
\sigma_j = z_{\rm b} \sqrt{k_{\rm w}^2 -\beta_j^2}, \quad \zeta_j =z_{\rm b} \sqrt{\beta_j^2-k_{\rm s}^2}  ,
\end{align}
and
\begin{align}
A_j^2 = \frac{2/z_{\rm b}}{ \frac{1}{\rho_{\rm w}} ( 1-\frac{\sin (2\sigma_j)}{2\sigma_j})
+\frac{1}{\rho_{\rm s}} \frac{\sin^2 (\sigma_j)}{\zeta_j} }  .
\end{align}
The $\sigma_j$'s are the solutions in $(0, \sqrt{k_{\rm w}^2-k_{\rm s}^2} z_{\rm b})$ of
\begin{align}
\tan (\sigma) = - \frac{\sigma}{\sqrt{ z_{\rm b}^2 (k_{\rm w}^2-k_{\rm s})^2-\sigma^2}} \frac{\rho_{\rm s}}{\rho_{\rm w}} ,
\end{align}
and we denote by $N$ the number of solutions.

{\bf Continuous spectrum.}
For $\gamma \in (-\infty,k_{\rm s}^2)$, the improper eigenvector has the form:
\begin{align}
\phi_\gamma(z) = 
\left\{
\begin{array}{l}
A_\gamma \sin (\eta_\gamma z/z_{\rm b}) \mbox{ if } 0\leq z\leq z_{\rm b}\\
A_\gamma\big[ \sin (\eta_\gamma ) \cos(\xi_\gamma (z-z_{\rm b})/z_{\rm b})
+\frac{\rho_{\rm s}}{\rho_{\rm w}} \frac{\eta_\gamma}{\xi_\gamma}\cos(\eta_\gamma) \sin(\xi_\gamma (z-z_{\rm b})/z_{\rm b})
\big]
\mbox{ if }  z\geq z_{\rm b}
\end{array}
\right.
\end{align}
where
\begin{align}
\eta_\gamma=z_{\rm b} \sqrt{k_{\rm w}^2-\gamma},\quad \xi_\gamma=z_{\rm b}\sqrt{k_{\rm s}^2-\gamma}
\end{align}
and
\begin{align}
A_\gamma^2 = \frac{\xi_\gamma \rho_{\rm s} z_{\rm b}}{\pi (\xi_\gamma^2 \sin^2(\eta_\gamma) +
\frac{\rho_{\rm s}^2}{\rho_{\rm w}^2} \eta_\gamma^2 \cos^2(\eta_\gamma))}  .
\end{align}
We remark that $\phi_\gamma$ does not belong to $L^2(\RR^+, \rho_0^{-1})$, but 
$\left( \phi_\gamma,\phi\right)$ can be defined for any test function $\phi \in L^2(\RR^+,\rho_0^{-1})$ as
\begin{align}
\left( \phi_\gamma,\phi\right) = \lim_{M \to +\infty} \int_0^M \phi_\gamma(z) \phi(z) \rho_0(z)^{-1} dz
,
\end{align}
where the limit holds on $L^2(-\infty,k_{\rm s}^2)$.

{\bf Completeness.}
We have for any $\phi\in L^2(\RR^+, \rho_0^{-1})$:
\begin{align}
\left( \phi,\phi \right)  = \sum_{j=1}^N \big| \left( \phi_j ,\phi\right) \big|^2
+\int_{-\infty}^{k_{\rm s}^2} \big| \left( \phi_\gamma,\phi\right)\big|^2 d\gamma   .
\end{align}
The map which assigns to every element of $L^2(\RR^+, \rho_0^{-1})$ the coefficients of its spectral decomposition
$$
\phi \mapsto \Big( \left( \phi_j,\phi\right)  ,j=1,\ldots,N, 
 \left(\phi_\gamma,\phi\right), \gamma\in ({-\infty},{k_{\rm s}^2})\Big)
$$
is an isometry from $L^2(\RR^+, \rho_0^{-1})$ onto $\CC^N \times L^2(-\infty,k_{\rm s}^2)$.
There exists a resolution of the identity $\Pi$ such that for any $\phi\in L^2(\RR^+ , \rho_0^{-1})$ and $-\infty \leq r \leq r' \leq +\infty$:
$$
\Pi(r,r')(\phi) = \sum_{j=1}^N \left( \phi_j ,\phi \right) \phi_j {\bf 1}_{(r,r')}(\beta_j^2)
+
\int_r^{\min(k_{\rm s}^2,r')} \left( \phi_\gamma, \phi  \right) \phi_\gamma d\gamma {\bf 1}_{(r,+\infty)}(k_{\rm s}^2) ,
$$
and for any $\phi$ in the domain of ${\cal H}$,
$$
\Pi(r,r')({\cal H}\phi)  = \sum_{j=1}^N \beta_j^2 \left( \phi_j ,\phi \right) \phi_j {\bf 1}_{(r,r')}(\beta_j^2)
+
\int_r^{\min(k_{\rm s}^2,r')} \gamma \left( \phi_\gamma, \phi \right) \phi_\gamma d\gamma {\bf 1}_{(r,+\infty)}(k_{\rm s}^2) .
$$

{\bf Modal decomposition.}
Any solution of the Helmholtz equation in homogeneous medium can be expanded as 
\begin{equation}
\label{eq:modalexpansion}
\hat{p}(x,z) = \sum_{j=1}^N \hat{p}_j(x) \phi_j(z)
+\int_0^{k_{\rm s}^2} \hat{p}_\gamma (x) \phi_\gamma(z) d\gamma
+\int_{-\infty}^0 \hat{p}_\gamma (x) \phi_\gamma(z) d\gamma .
\end{equation}
The first modes are guided, the second ones  are radiating, the third ones are evanescent.
Indeed, the complex mode amplitudes satisfy
\begin{align}
\partial_x^2 \hat{p}_j +\beta_j^2 \hat{p}_j&=0, \quad  j=1,\ldots,N,\\
\partial_x^2 \hat{p}_\gamma +\gamma \hat{p}_\gamma&=0, \quad  \gamma \in (-\infty,k_{\rm s}^2) .
\end{align}
Therefore, if the source is in the plane $x=0$, at $(0,z_0)$, as in (\ref{eq:pressure}), we have for $x>0$:
\begin{align}
\nonumber
\hat{p}(x,z) =&  \sum_{j=1}^N \frac{\hat{a}_{j,0}}{\sqrt{\beta_j}} e^{i \beta_j x} \phi_j(z)
+\int_0^{k_{\rm s}^2}  \frac{\hat{a}_{\gamma,0}}{\gamma^{1/4}} e^{i \sqrt{\gamma} x}\phi_\gamma(z) d\gamma
\\
&+\int_{-\infty}^0  \frac{\hat{a}_{\gamma,0}}{|\gamma|^{1/4}}   e^{- \sqrt{|\gamma|} x} \phi_\gamma(z) d\gamma
,
\end{align}
with the modal amplitudes determined by the source:
\begin{align}
\hat{a}_{j,0} =  & \frac{\sqrt{\beta_j}}{2} \phi_j(z_0)  ,\quad j=1,\ldots,N,\label{a_j0}\\
\hat{a}_{\gamma,0} = & \frac{\gamma^{1/4}}{2}   \phi_\gamma(z_0),\quad \gamma \in (0,k_{\rm s}^2) ,\\
\hat{a}_{\gamma,0} = & \frac{|\gamma|^{1/4}}{2}  \phi_\gamma(z_0) ,\quad \gamma \in (-\infty,0) .
\end{align}

\section{Random and dissipative waveguide}
\label{sec:randomwaveguide}
Here we assume that the waveguide is weakly randomly perturbed and weakly dissipative:
$$
c^2(x,z) = \frac{c_0^2(z)}{1+V(x,z)} .
$$
The perturbation $V(x,z)$ has two components: random velocity perturbation in water
 and deterministic constant dissipation in sediments:
\begin{align}
V(x,z) =  \nu(x,z) {\bf 1}_{(0,z_{\rm b})}(z)
+ i  \nu_{\rm s}  {\bf 1}_{ (z_{\rm b},+\infty)}(z),
\end{align}
where $\nu(x,z)$ is a zero-mean random process that describes the relative fluctuations
of the propagation speed in water and $\nu_{\rm s}>0$ models the damping in the sediments.

\subsection{The coupled mode equations}
The solution of the perturbed Helmholtz equation
\begin{align}
\Big[(\partial_z^2 + \partial_x^2 )+ \frac{\omega^2}{c_0^2(z)} + \frac{\omega^2}{c_0^2(z)} V(x,z)\Big] \hat{p}( x,z) =  
\delta(x)\delta(z-z_0) ,  
\end{align}
for $x\in \RR, z\in (0,z_{\rm b})\cup(z_{\rm b},+\infty)$, 
with the boundary conditions (\ref{eq:bc}),
can be expanded as (\ref{eq:modalexpansion}) and the complex mode amplitudes satisfy the coupled equations:
\begin{align}
\partial_x^2 \hat{p}_j +\beta_j^2 \hat{p}_j&= -\omega^2 \sum_{l=1}^N C_{jl}(x) \hat{p}_l 
-\omega^2 \int_{-\infty}^{k_{\rm s}^2} C_{j\gamma'} (x) \hat{p}_{\gamma'} d\gamma', \quad  j=1,\ldots,N,\\
\partial_x^2 \hat{p}_\gamma +\gamma \hat{p}_\gamma&= -\omega^2 \sum_{l=1}^N C_{\gamma l}(x) \hat{p}_l 
-\omega^2 \int_{-\infty}^{k_{\rm s}^2} C_{\gamma\gamma'} (x) \hat{p}_{\gamma'} d\gamma',  
\quad  \gamma \in (-\infty,k_{\rm s}^2)  ,
\end{align}
with
\begin{align}
C_{jl}(x) =& \left(\phi_j,\phi_l V(x,\cdot) c_0^{-2}\right) ,\\
C_{j\gamma'}(x) =& \left(\phi_j,\phi_{\gamma'} V(x,\cdot) c_0^{-2}  \right) ,\\
C_{\gamma l}(x) =&  \left(\phi_\gamma,\phi_l V(x,\cdot)c_0^{-2} \right) ,\\
C_{\gamma\gamma'}(x) =& \left(\phi_\gamma,\phi_{\gamma'} V(x,\cdot) c_0^{-2}  \right) .
\end{align}
The coupling term $C_{jl}$ has two components:
\begin{align*}
C_{jl}(x) =  C^{\rm w}_{jl}(x) + i  \nu_{\rm s} \,\rho_{\rm s}^{-1} c_{\rm s}^{-2} \int_{z_{\rm b}}^{+\infty}  \phi_j(z) \phi_l(z) dz ,
\end{align*}
where
\begin{align*}
 C^{\rm w}_{jl}(x) = c_{\rm w}^{-2} \rho_{\rm w}^{-1} \int_0^{z_{\rm b}}  \nu(x,z) \phi_j(z) \phi_l(z) dz ,
\end{align*}
and similarly for the other coupling terms.

We introduce the amplitudes of the generalized right and left-going mode amplitudes $\hat{a}_j$ and $\hat{b}_j$
for the guided and radiating modes:
\begin{align}
\hat{p}_j(x) =& \frac{1}{\sqrt{\beta_j}}\Big( \hat{a}_j (x) e^{i\beta_j x} +\hat{b}_j(x) e^{- i\beta_j x} \Big), \quad j=1,\ldots,N,\\
\partial_x  \hat{p}_j(x) =& i\sqrt{\beta_j}\Big( \hat{a}_j (x) e^{i\beta_j x} -\hat{b}_j(x) e^{- i\beta_j x} \Big),\quad j=1,\ldots,N,\\
\hat{p}_\gamma(x) =& \frac{1}{\gamma^{1/4}}\Big( \hat{a}_\gamma (x) e^{i\sqrt{\gamma} x} 
+\hat{b}_\gamma(x) e^{- i\sqrt{\gamma}x} \Big),\quad \gamma\in (0,k_{\rm s}^2) ,\\
\partial_x  \hat{p}_\gamma(x) =& i\gamma^{1/4}\Big( \hat{a}_\gamma (x) e^{i\sqrt{\gamma} x} 
-\hat{b}_\gamma(x) e^{- i\sqrt{\gamma}x} \Big),\quad \gamma\in (0,k_{\rm s}^2) .
\end{align}
They satisfy the following coupled differential equations:
\begin{align*}
\frac{d}{dx} \hat a_j(x) =& \frac{i \omega^2}{2} \sum_{l =1}^N \frac{C_{jl}(x)}{\sqrt{\beta_j \beta_l}}(\hat a_l( x) e^{i(\beta_l-\beta_j) x} + \hat b_l( x) e^{-i (\beta_l+\beta_j) x}) \\
&+ \frac{i \omega^2}{2} \int_{0}^{k_{\rm s}^2} \frac{C_{j \gamma'}(x)}{\sqrt{\beta_j \sqrt{\gamma'}}}(\hat a_{\gamma'}( x) e^{i(\sqrt{\gamma'}-\beta_j) x} + \hat b_{\gamma'}(x) e^{-i (\sqrt{\gamma'}+\beta_j) x}) d\gamma' \\
&+\frac{i \omega^2}{2} \int_{-\infty}^{0} \frac{C_{j \gamma'}(x)}{\sqrt{\beta_j}}\hat p_{\gamma'}( x) e^{- i \beta_j x} d\gamma'  ,
\end{align*}
\begin{align*}
\frac{d}{dx} \hat b_j( x) =& - \frac{i \omega^2}{2} \sum_{l =1}^N \frac{C_{jl}(x)}{\sqrt{\beta_j \beta_l}}(\hat a_l( x) e^{i(\beta_l+\beta_j) x} + \hat b_l(x) e^{-i (\beta_l-\beta_j) x}) \\
&- \frac{i \omega^2}{2} \int_{0}^{k_{\rm s}^2} \frac{C_{j \gamma'}(x)}{\sqrt{\beta_j \sqrt{\gamma'}}}(\hat a_{\gamma'}(x) e^{i(\sqrt{\gamma'}+\beta_j) x} + \hat b_{\gamma'}( x) e^{- i (\sqrt{\gamma'}-\beta_j) x}) d\gamma' \\
&-\frac{i \omega^2}{2} \int_{-\infty}^{0} \frac{C_{j \gamma'}(x)}{\sqrt{\beta_j}}\hat p_{\gamma'}(x) e^{i \beta_j x} d\gamma'  .
\end{align*}
The radiating  mode amplitudes $\hat{a}_\gamma (x)$ and $\hat{b}_\gamma (x)$, $\gamma \in (0,k_{\rm s}^2)$,
satisfy the same equations upon replacing $j$ by $\gamma$ in the formulas above.
The evanescent mode amplitudes $\hat p_{\gamma}(x)$, $\gamma \in (-\infty,0)$, satisfy the differential equations
\begin{align}
\label{eq:evan1}
\partial_x^2 \hat{p}_{\gamma} +\gamma \hat{p}_{\gamma} = 
-  g_{\gamma}(x) -\hat{g}_{\gamma}^{\rm ev}(x) , \quad x \ne 0, 
\end{align}
where $\hat{g}_{\gamma}(x)$ and  $\hat{g}_{\gamma}^{\rm ev}(x)$  
are given by
\begin{align}
\nonumber
\hat{g}_{\gamma}(x) =& \omega^2
\sum_{l'=1}^{N} \frac{C_{\gamma l'}(x)}{\sqrt{\beta_{l'}}} \Big[ {a}_{l'}(x) e^{i \beta_{l'} x}
+ {b}_{l'}(x) e^{- i\beta_{l'} x}
\Big]\\
&+   \omega^2 
\int_0^{k_{\rm s}^2}  \frac{C_{\gamma\gamma'}(x)}{\sqrt[4]{\gamma'}} \Big[ {a}_{\gamma'}(z) e^{i  \sqrt{\gamma'} x}
+ {b}_{\gamma'}(x) e^{- i \sqrt{\gamma'} z}
\Big]d\gamma',  \label{eq:g_tg}\\
\hat{g}^{\rm ev}_{\gamma}(x) =&
 \omega^2  \int_{-\infty}^0 C_{\gamma\gamma'}(x) \hat{p}_{\gamma'}  (x)
d\gamma' .
\label{eq:defgtga}
\end{align}

We wish to understand the evolution of the complex mode amplitudes in the regime where the perturbations are weak and 
the propagation distance is large. We assume 
\begin{align}
V(x,z) =  \eps \nu(x,z) {\bf 1}_{z\in (0,z_{\rm b})}
+ i \eps^2 \nu_{\rm s}  {\bf 1}_{z\in (z_{\rm b},+\infty)},
\end{align}
where $\eps$ is a small dimensionless parameter that expresses the fact that the random perturbations are weak and dissipation is even weaker.
We suppose that we have a perfect (homogeneous and non-dissipative) waveguide outside the region $x\in(0,x_\infty/\eps^2)$ for some $x_\infty>0$.
We carry out an asymptotic study in the small $\eps$ limit. 
In this framework the evanescent mode amplitudes can be expressed in terms of the propagating mode amplitudes \cite{garnier_evan,gomez}
and the differential equation for $\hat{a}_j^{\eps}(\omega,x):=\hat{a}_j(\omega,x/\eps^2)$ can be written as:
\begin{align*}
\frac{d}{dx} \hat{a}_j^{\eps}(x) &= \frac{i \omega^2}{2} \sum_{l =1}^N \frac{\eps^{-1} C^{\rm w}_{jl}(x/\eps^2) + i \nu_{\rm s} \rho_{\rm s}^{-1} c_{\rm s}^{-2} \int_{z_{\rm b}}^{\infty} \phi_j(z) \phi_l(z) dz}{\sqrt{\beta_j \beta_l}} \\
&\quad\quad \times  (\hat{a}_l^{\eps}(x) e^{i(\beta_l-\beta_j) x/\eps^2}+   \hat{b}_l^\eps(x) e^{-i (\beta_l+\beta_j) x/\eps^2}) + O(1).
\end{align*}
In the asymptotic framework ``weak fluctuations, weak dissipation, large propagation distance", 
we carry out a multiscale analysis similar as the one reported in Ref.~~\onlinecite{gomez}.
We find that the random fluctuations of the propagation speed in water give the standard coupling terms
between the right-going mode amplitudes of the guided modes. They also induce loss towards the radiating modes.
Moreover, the dissipation in the sediments gives rise to additional damping terms in the evolution equations
of the complex mode amplitudes.
The coupled system of equations for the complex mode amplitudes is of the form $d{\itbf a}^\eps/dx = \eps^{-1} {\itbf F}(x/\eps^2, {\itbf a}^\eps) + {\itbf G}(x/\eps^2, {\itbf a}^\eps)$. 
The small $\eps$ limit of such differential equations are studied in the book~~\onlinecite{FGPSbook} (Chapters 6 and 20). 
The full result for the complex mode amplitudes is expressed in Appendix \ref{app:A}.
From this result we get that the process $(|\hat{a}_j^\eps(x)|^2)_{j=1}^N$ converges towards a Markov process ${\itbf P}(x)=(P_j(x))_{j=1}^N$ whose infinitesimal generator $\mathcal{L}_{\itbf P}$ writes:
\begin{align}
\label{eq:genLP}
\mathcal{L}_{\itbf P} = \sum_{j \neq l}  \Gamma_{lj} \Big[P_l P_j (\frac{\partial}{\partial P_j} - \frac{\partial}{\partial P_l}) \frac{\partial}{\partial P_j} + (P_l - P_j) \frac{\partial}{\partial P_j}\Big] -  \sum_{j=1}^N \Lambda_j P_j \frac{\partial}{\partial P_j} ,
\end{align}
with
\begin{align}
\nonumber
\label{def:Lambdaj}
\Lambda_j =& \int_0^{k_{\rm s}^2} \frac{\omega^4}{2\sqrt{\gamma}\beta_j}
\int_0^\infty \EE [ C^{\rm w}_{j\gamma}(0)C^{\rm w}_{j\gamma}(x) ] \cos \big((\sqrt{\gamma}-\beta_j)x\big) dx d\gamma
\\
&+
\nu_{\rm s} \frac{\omega^2}{\beta_j c_{\rm s}^2 \rho_{\rm s}} \int_{z_{\rm b}}^\infty \phi_j(z)^2 dz  , \\
\Gamma_{lj}=& \frac{\omega^4}{2 \beta_j \beta_l}
\int_0^\infty \EE [C^{\rm w}_{jl}(0)C^{\rm w}_{jl}(x)] \cos \big((\beta_l-\beta_j)x\big) dx  \quad \mbox{ for }j \neq l,
\label{def:Gammalj}
\end{align}
where
\begin{align}
\EE [ C^{\rm w}_{j\gamma}(0)C^{\rm w}_{j\gamma}(x) ]
=& 
\rho_{\rm w}^{-2}  c_{\rm w}^{-4}
\int_0^{z_{\rm b}}\int_0^{z_{\rm b}} \phi_j \phi_\gamma(z) \EE [\nu(0,z)\nu(x,z')] \phi_j \phi_\gamma(z') dz dz' ,\\
\EE [C^{\rm w}_{jl}(0)C^{\rm w}_{jl}(x)]  =&
\rho_{\rm w}^{-2}  c_{\rm w}^{-4}
\int_0^{z_{\rm b}}\int_0^{z_{\rm b}} \phi_j \phi_l(z) \EE [\nu(0,z)\nu(x,z')] \phi_j \phi_l(z') dz dz' ,
\end{align}
and we have $ \int_{z_{\rm b}}^\infty \phi_j(z)^2 dz = A_j^2 \sin^2(\sigma_j) \,z_{\rm b}/(2\zeta_j)$.
This means that, for any test function $F:\mathbb{R}^N \to \mathbb{R}$,
the expectation $\EE[ F( {\itbf P}(x) )]$ satisfies the Kolmogorov forward equation
\begin{equation}
\label{eq:fk1}
\partial_x \EE [ F( {\itbf P}(x) )] = \EE [{\cal L}_{\itbf P}F( {\itbf P}(x) )] .
\end{equation}
Equivalently, the probability density function $p(x,{\itbf P})$ of the random vector ${\itbf P}(x)$ satisfies the Kolmogorov backward
(or Fokker-Planck) equation
$
\partial_x p=  {\cal L}_{\itbf P}^* p ,
$
where ${\cal L}_{{\itbf P}}^*$ is the adjoint of $ {\cal L}_{\itbf P}$.

From the form of the generator ${\cal L}_{\itbf P}$, one can establish that the $n$th-order moments 
of the mode powers satisfy closed equations.
In particular, using (\ref{eq:genLP}) and (\ref{eq:fk1}) we find that the mean mode powers
\begin{equation}
{Q}_j(x)=
\EE [P_j(x)]  = \lim_{\eps \to 0} \EE \big[ |\hat{a}_j^\eps(x)|^2\big]
\end{equation}
satisfy the closed system of equations
\begin{equation}
\partial_x {Q}_j = -\Lambda_j {Q}_j +\sum_{l=1}^N \Gamma_{lj} \big({Q}_l-{Q}_j\big)  ,
\quad {Q}_j(0)=|\hat{a}_{j,0}|^2 .
\end{equation}
The form of these coupled-mode equations is well-known \cite{dozier} although the mode-dependent attenuation term $\Lambda_j$ 
was usually introduced heuristically so far.
The solution explicitly writes:
\begin{align}
\label{eq:mom1b}
{\itbf Q}(x) = \exp({\bf A} x) {\itbf Q}(0)  ,
\end{align}
with the matrix ${\bf A}$ defined by ($\delta_{jl}$ is the Kronecker symbol and $\Gamma_{jj} = -\sum_{l'\neq j}\Gamma_{jl'} $):
\begin{align*}
{\bf A} := (\Gamma_{jl} - \Lambda_j \delta_{jl})_{j,l=1}^N .
\end{align*}

\subsection{Computation of the coefficients of the matrix ${\bf A}$}
The goal of this paragraph is to show how to get closed-form expressions of the coefficients of
 the matrix ${\bf A}$ when the correlation of the perturbation decreases exponentially 
as a function of the horizontal distance between the two points, i.e.
\begin{align*}
\EE [\nu(x,z) \nu(x',z')] = \sigma^2 \exp(-|x-x'|/\ell_{\rm h}) R(z,z') ,
\end{align*}
where $\ell_{\rm h}$ is the horizontal correlation radius of the random fluctuations of the index of refraction.
We can compute $\Gamma_{jl}$ for $j \neq l$ from:
\begin{align*}
&\int_0^\infty \EE [C^{\rm w}_{jl}(0)C^{\rm w}_{jl}(x)] \cos \big((\beta_l-\beta_j)x\big) dx \\
&\quad = \frac{\sigma^2}{\rho_{\rm w}^{2} c_{\rm w}^{4}} \frac{\ell_{\rm h}}{1+ (\beta_l - \beta_j)^2 \ell_{\rm h}^2}  \int_0^{z_{\rm b}} \int_0^{z_{\rm b}} R(z,z') \phi_j \phi_l(z) \phi_j \phi_l(z') dz dz'
\end{align*}
and similarly we can compute $\Lambda_j$ from the same quantities upon substitution $\gamma$ for $l$.
By (\ref{def:Gammalj}) we have
\begin{align*}
\Gamma_{lj}=& \frac{\omega^4}{2 \beta_j \beta_l} \frac{\sigma^2}{\rho_{\rm w}^{2} c_{\rm w}^{4}} \frac{\ell_{\rm h}}{1+ (\beta_l - \beta_j)^2 \ell_{\rm h}^2} \int_0^{z_{\rm b}} \int_0^{z_{\rm b}} R(z,z') \phi_j \phi_l(z) \phi_j \phi_l(z') dz dz' .
\end{align*}
Replacing $\phi_j$ for all $j$ by their expression, and using the notation $k_{{\rm w}j} := \sqrt{k_{\rm w}^2 -\beta_j^2}$, we obtain: 
\begin{align}
\nonumber
 &\int_0^{z_{\rm b}} \int_0^{z_{\rm b}} R(z,z') \phi_j \phi_l(z) \phi_j \phi_l(z') dz dz' 
 = \frac{A_j^2 A_l^2}{4} \big(S(k_{{\rm w}j} - k_{{\rm w}l}, k_{{\rm w}j} - k_{{\rm w}l}) \\
 & + S(k_{{\rm w}j} + k_{{\rm w}l}, k_{{\rm w}j} + k_{{\rm w}l})  - S(k_{{\rm w}j} - k_{{\rm w}l}, k_{{\rm w}j} + k_{{\rm w}l}) - S(k_{{\rm w}j} + k_{{\rm w}l}, k_{{\rm w}j} - k_{{\rm w}l})\big),
  \label{eq:Gjl}
\end{align}
where 
\begin{align*}
S(k, k') :=  \int_0^{z_{\rm b}}\int_0^{z_{\rm b}} R(z,z') \cos(k\, z) \cos(k'\, z') dz dz'.
\end{align*}
When $R$ has the exponential form $R(z,z') = \exp(-|z-z'|/\ell_{\rm v})/2$, then
the variance of the fluctuations is $\sigma^2/2$ and $S(k,k')$
has a closed-form expression. 
We will use this particular form of the correlation function of the medium in Section \ref{sec:ip} because it gives simple expressions
for the coefficients $\Gamma_{jl}$ and $\Lambda_j$, which is convenient for the resolution of the inverse problem
which requires many evaluations of such coefficients for different values $\sigma, \ell_{\rm v},\ell_{\rm h}$ of the statistics of the random medium 
and $\rho_{\rm s}, \nu_{\rm s}$ of the sea bottom.
We could as well use the  Garrett-Munk correlation function at the expense of some computational overburden \cite{garrett},
but the impact of the exact form of the correlation function turns out to be weak.

\subsection{Pressure field correlations}
\label{subsec:verticalcorrelations}
We now compute the correlation of the received signal on a vertical segment at a fixed horizontal distance $x_{\rm a}$. When $x_{\rm a} \gg 1$, we can use the asymptotic study from above and the decomposition of $\hat p( x_{\rm a}, z)$ in the basis of eigenvectors $\phi_j$ is: 
\begin{align*}
\hat p( x_{\rm a}, z) \simeq \sum_{j=1}^N \frac{1}{\sqrt{\beta_j}} \hat{a}_j (x_{\rm a})  \phi_j(z).
\end{align*}
Therefore
\begin{align*}
\EE \big[\overline{ \hat{p}( x_{\rm a}, z)}  \hat{p}( x_{\rm a}, z') \big] 
&= \sum_{j=1}^N \frac{\EE [|\hat a_j(x_{\rm a})|^2]}{\beta_j} \phi_j(z)\phi_j(z')+\sum_{j\neq l} \frac{\EE [\overline{\hat a_j(x_{\rm a})}\hat a_l(x_{\rm a}) ]}{\sqrt{\beta_j \beta_l}} \phi_j(z)\phi_l(z').
\end{align*}
The cross second moments $\EE [\overline{\hat a_j }\hat a_l  ]$ for $j\neq l$ decay exponentially with the propagation distance $x$
as shown in Ref.~~\onlinecite{GP07} (see also Ref.~~\onlinecite{colosi12}) and they can be neglected as soon as the propagation distance $x$ is larger than the scattering mean free path,
that depends on the coefficients $\Gamma_{jl}$. This is roughly speaking the distance beyond which the phase of the 
field has a random part which variance of the order of or larger than one, so that the coherent field is vanishing \cite{GP07,borcea15}.
The second moments $\EE [|\hat a_j |^2 ]$ were computed
 in the previous paragraph using the exponential of the matrix ${\bf A}$:
\begin{align*}
\EE [|\hat a_j(x_{\rm a})|^2] =  \sum_{l=1}^N (\exp({\bf A} x_{\rm a}))_{jl} |\hat a_{l,0}|^2,
\end{align*}
where $\hat a_{l,0}$ is given by $\eqref{a_j0}$.
We can now write the correlation between the sound pressure signals recorded by two receivers which are separated by $y$ along the vertical segment at distance $x_{\rm a}$ and of depth $[z_m,z_M]$ where $0 <z_m < z_M < z_{\rm b}$.
For all $y \in [0,z_M -z_m]$, the spatial correlation $C_{x_{\rm a}}(y)$ writes:
\begin{align}
\nonumber
C_{x_{\rm a}}(y) &:= \frac{1}{z_M-z_m -y} \int_{z_m}^{z_M-y} \EE \big[\overline{ \hat{p}( x_{\rm a}, z)} \hat{p}( x_{\rm a}, z+y) \big] dz \\
\nonumber
&= \frac{1}{z_M-z_m -y} \sum_j \frac{\EE [|\hat a_j(x_{\rm a})|^2]}{\beta_j} \frac{A_j^2}{2} 
 \Big(\cos(k_{{\rm w}j} y) (z_M-z_m-y) \\
                 & \quad - \frac{\sin(k_{{\rm w}j}(2 z_M-y)) - \sin(k_{{\rm w}j}(2 z_m+y))}{2 k_{{\rm w}j}}\Big).
                   \label{cst}
\end{align}

\subsection{Equipartition regime}
By (\ref{eq:mom1b}) the mean mode powers satisfy
$$
{Q}_j(x) \stackrel{x \to +\infty}{\simeq}
c_V V_j  \exp\big( -\lambda  x \big) \big(1+o(1)\big)  ,
$$
where $({\itbf V} ,-\lambda )$ is the first eigenvector/eigenvalue of the matrix
${\bf A}=\boldsymbol{\Gamma} -\boldsymbol{\Phi} $,
with $\Gamma_{jl}$ given by (\ref{def:Gammalj}) for $j\neq l$,  $\Gamma_{jj} = -\sum_{l'\neq j}\Gamma_{jl'} $, 
$\Phi_{jl} = \Lambda_j \delta_{jl}$,
and
$$
c_V = \sum_{j=1}^N V_j  |\hat{a}_{j,0}|^2  .
$$
Note that:
\begin{itemize}
\item The coefficients of ${\itbf V}$ have all the same sign (so we can assume that they are nonnegative).\\
Indeed, the solution of the equation $\partial_x {\itbf Q} = (\boldsymbol{\Gamma}-\boldsymbol{\Phi}) {\itbf Q}$ starting from any
${\itbf Q}(x=0)$ with nonnegative coefficients must have nonnegative coefficients since the exponential of a matrix whose off diagonal entries are nonnegative has nonnegative coefficients.
As ${\itbf Q}(x)$ is equivalent to $({\itbf V}^T{\itbf Q}(x=0)) 
\exp(-\lambda  x) {\itbf V}$, the coefficients of ${\itbf V}$ must have the same sign. We deduce as well that $c_V$ has the same sign as ${\itbf V}$.
\item $\lambda \geq0$.

Indeed, we have $ (\boldsymbol{\Gamma} -\boldsymbol{\Phi} )  {\itbf V} = - \lambda {\itbf V} $. By projecting onto
${\itbf V}^{(0)} = \big(1/\sqrt{N} \big)_{j=1}^N$, which is in the kernel of $\boldsymbol{\Gamma} $, we get
${{\itbf V}^{(0)}}^T \boldsymbol{\Phi}   {\itbf V} = \lambda   {{\itbf V}^{(0)}}^T {\itbf V}$.
Since $\boldsymbol{\Phi} $ is diagonal with nonnegative coefficients and the coefficients of ${\itbf V} $ have the same sign and
cannot be all zero, we find $ \lambda  = \sum_j \Lambda_j {V}_j / \sum_j {V}_j \geq 0$.
\end{itemize}

In the following we discuss cases with zero or weak dissipation, where explicit expressions can be obtained.
Remember that the effective dissipation is the sum of the two effects: power leakage from the guided modes to the radiating modes and
dissipation in the sediments.

{\bf No effective dissipation.}
If there is no effective dissipation $\boldsymbol{\Phi} ={\bf 0}$, 
then the first eigenvector/eigenvalue $({\itbf V}^{(0)},-\lambda^{(0)})$ of the matrix
$\boldsymbol{\Gamma} $ is 
$$
{\itbf V}^{(0)} = \big(1/\sqrt{N} \big)_{j=1}^N, \quad
\lambda^{(0)}=0 ,
$$
which gives the standard equipartition result \cite{creamer,FGPSbook,GP07}:
$$
{Q}_j(x) \stackrel{x \to +\infty}{\longrightarrow}
\frac{1}{N}  \sum_{j=1}^N  |\hat{a}_{j,0}|^2   .
$$

{\bf Weak effective dissipation.}
We next consider the case when the effective dissipation is weak, 
that is to say,  the matrix $\boldsymbol{\Phi} $ is much smaller
  than the matrix $\boldsymbol{\Gamma} $, with a typical ratio of the order of $\delta$.
We then assume that
$\Lambda_j = \delta \Lambda_j^{(1)}  ,$
with $\delta \ll 1$.
Then we can write $\boldsymbol{\Phi}  =\delta \boldsymbol{\Phi}^{(1)}$ 
with $\Phi_{jl}^{(1)} =
\Lambda_j^{(1)} \delta_{jl}$
and $\boldsymbol{\Gamma}  = \boldsymbol{\Gamma}^{(0)}$
and the first eigenvector/eigenvalue $({\itbf V} ,-\lambda)$  of the matrix
$\boldsymbol{\Gamma} -\boldsymbol{\Phi} $ can be expanded
as
$$
 {\itbf V}  = {\itbf V}^{(0)} +\delta {\itbf V}^{(1)} +O(\delta^2)
  ,\quad  \lambda= \delta \lambda^{(1)} +\delta^2 \lambda^{(2)} +O(\delta^3),
$$
with
\begin{eqnarray}
\lambda^{(1)} &=& {{\itbf V}^{(0)}}^T \boldsymbol{\Phi}^{(1)} {\itbf V}^{(0)}
=
\frac{1}{N}\sum_{j=1}^N \Lambda_j^{(1)},\\
\lambda^{(2)} &=& {{\itbf V}^{(0)}}^T \boldsymbol{\Gamma}^{(0)} {\itbf V}^{(1)} ,
\end{eqnarray}
and $ {\itbf V}^{(1)}$ is solution of $ \boldsymbol{\Gamma}^{(0)}  {\itbf V}^{(1)} = (\boldsymbol{\Phi}^{(1)} - \lambda^{(1)} )
{\itbf V}^{(0)}$ and is orthogonal to $ {\itbf V}^{(0)}$.
If, for instance, $\Gamma_{jl}\equiv \Gamma >0$ for all $j\neq l$, then
$$
{\itbf V}^{(1)} = -\frac{1}{\Gamma N^{3/2}} \big(\Lambda_j^{(1)} \big)_{j=1}^N
$$
and
$$
\lambda^{(2)} = 
- \frac{1}{\Gamma N^2} \sum_{j=1}^N (\Lambda_j^{(1)} -\lambda^{(1)})^2 .
$$

\subsection{Fluctuation analysis}
By (\ref{eq:genLP}) we find that the second-order moments of the mode powers 
$$
R_{jl}  (x)= \EE \big[P_j(x)P_l(x)\big] ,\quad j,l=1,\ldots,N,
$$
satisfy the closed equations
\begin{align}
&
 \partial_x R_{jj} = - 2 \Lambda_j R_{jj}  +\sum_{n \neq j}  \Gamma_{jn} (4R_{jn} -2R_{jj} )  ,\\
&
 \partial_x R_{jl}  = - (2\Gamma_{jl}+ \Lambda_j +\Lambda_l) R_{jl}  +
\sum_{n \neq l}  \Gamma_{ln} (R_{jn} -R_{jl} )  + \sum_{n \neq j}  \Gamma_{jn} (R_{nl} -R_{jl} )  ,
\end{align}
for $j\neq l$. 
This system has the same form as the one found in the literature dedicated to coupled mode theory \cite{dozier,creamer}.
The initial conditions are $ R_{jl} (0)=|\hat{a}_{j,0}|^2|\hat{a}_{l,0}|^2 $.
Let us introduce ${\itbf S} = (S_{jl})_{1\leq j\leq l\leq N}$ defined by
\begin{equation}
S_{jl} = 
\left\{
\begin{array}{ll}
R_{jl}+R_{lj} &\mbox{ if } j  <  l,\\
R_{jj} & \mbox{ if } j=l  .
\end{array}
\right.
\end{equation}
The $S_{jl}$'s satisfy the system
\begin{align}
\nonumber
\partial_x S_{jl} =& - (\Lambda_j+\Lambda_l)S_{jl}
+\sum_{n \not\in \{j,l\}} \big[ \Gamma_{ln} (S_{jn}-S_{jl})+\Gamma_{jn}(S_{nl}-S_{jl})\big] \\
&
+2\Gamma_{jl} {\bf 1}_{j \neq l} (S_{jj}+S_{ll} -2S_{jl}) ,
\end{align}
with the convention that whenever $S_{jl}$ occurs with $j>l$, it is replaced by $S_{lj}$.
This can be written in the form $\partial_x {\itbf S}  = (\boldsymbol{\Theta} 
- \boldsymbol{\Psi} ) {\itbf S}$.
The linear operator $ \boldsymbol{\Psi} $ (that depends on the $\Lambda_j$'s) is diagonal and the linear operator $\boldsymbol{\Theta} $
(that depends on the $\Gamma_{jl}$'s)  is self-adjoint:  for any ${\itbf T}$ and $\widetilde{\itbf T}$, we have
\begin{align*}
\sum_{j\leq l}  (\boldsymbol{\Theta}  {\itbf T})_{jl}  \widetilde{T}_{jl}
&= -\sum_{j\leq l} \Theta_{jl,jl} T_{jl} \widetilde{T}_{jl}
+
\sum_{j <l , n \not\in \{j,l\}} \big( \Gamma_{ln} T_{jn} \widetilde{T}_{jl} +\Gamma_{jn}T_{nl}\widetilde{T}_{jl}\big) \\
&\quad +\sum_{j \neq n}\big(\Gamma_{jn} T_{jn} \widetilde{T}_{jj} +\Gamma_{jn} T_{nj} \widetilde{T}_{jj}\big) 
+2 \sum_{j < l}\big( \Gamma_{jl} T_{jj} \widetilde{T}_{jl} + \Gamma_{jl} T_{jj} \widetilde{T}_{jl}\big)\\
&=\sum_{j\leq l}
T_{jl} (\boldsymbol{\Theta}  \widetilde{\itbf T})_{jl} 
,
\end{align*}
because $2\sum_{j<l} = \sum_{j \neq l}$.
As a consequence, $\boldsymbol{\Theta} 
- \boldsymbol{\Psi} $ can be diagonalized and 
we find that
$$
{\itbf S}(x) \stackrel{x \to +\infty}{\simeq}
c_W {\itbf W} \exp\big( -\mu x \big) \big(1+o(1)\big)  ,
$$
where $c_W$ is the projection on the first eigenvector $ {\itbf W}$ in the basis of eigenvectors of $\boldsymbol{\Theta}-\boldsymbol{\Psi}$,
and
$$
c_W = 
2 \sum_{j< l} W_{jl}  |\hat{a}_{j,0}|^2|\hat{a}_{l,0}|^2  +  \sum_{j=1}^N W_{jj}  |\hat{a}_{j,0}|^4
=
 \sum_{j,l=1}^N W_{jl}  |\hat{a}_{j,0}|^2|\hat{a}_{l,0}|^2  ,
$$
with the convention that whenever $W_{jl}$ occurs with $j>l$, it is replaced by $W_{lj}$.

{\bf No effective dissipation.}
If there is no effective dissipation, then the first eigenvector/eigenvalue $({\itbf W}^{(0)},-\mu^{(0)})$ of the matrix
$\boldsymbol{\Theta}$ is 
$$
{\itbf W}^{(0)} = \big( c_N \big)_{1\leq j \leq l \leq N}, \quad
\mu^{(0)}=0 ,
$$
with $c_N=\sqrt{2}/\sqrt{N(N+1)}$ 
and 
$$
{\itbf S}(x) \stackrel{x \to +\infty}{\to}
c_W {\itbf W^{(0)}}  .
$$
As $\sum_{j\leq l} S_{jl}(x) =\sum_{j,l} R_{jl}(x) = (\sum_{j=1}^N |\hat{a}_{j,0}|^2)^2$, we deduce 
$$
{\itbf S}(x) \stackrel{x \to +\infty}{\to} \Big(\sum_{j=1}^N |\hat{a}_{j,0}|^2\Big)^2\frac{2}{N(N+1)} ,
$$
and
$$
{\itbf R}(x) \stackrel{x \to +\infty}{\to} \Big(\sum_{j=1}^N |\hat{a}_{j,0}|^2\Big)^2\frac{1+\delta_{jl}}{N(N+1)}.
$$

{\bf Weak effective dissipation.}
We next consider the case when the effective dissipation is weak, 
say $\Lambda_j = \delta \Lambda_j^{(1)}$.
Then we can write $\boldsymbol{\Psi} =\delta \boldsymbol{\Psi}^{(1)}$ and $\boldsymbol{\Theta} = \boldsymbol{\Theta}^{(0)}$
and the first eigenvector/eigenvalue $({\itbf W},-\mu)$  of the matrix
$\boldsymbol{\Theta}-\boldsymbol{\Psi}$ can be expanded
as
$$
 {\itbf W} = {\itbf W}^{(0)} +\delta {\itbf W}^{(1)} +O(\delta^2)
  ,\quad  \mu= \delta \mu^{(1)} +\delta^2 \mu^{(2)} +O(\delta^3),
$$
with
\begin{eqnarray}
\nonumber
\mu^{(1)} &=& {{\itbf W}^{(0)}}^T \boldsymbol{\Psi}^{(1)} {\itbf W}^{(0)}
=
c_N^2 \sum_{j <l} (\Lambda_j^{(1)}+\Lambda_l^{(1)})
+
c_N^2  \sum_{j} 2 \Lambda_j^{(1)}\\
 &=& \frac{2}{N} \sum_{j=1}^N  \Lambda_j^{(1)} = 2\lambda^{(1)},\\
\mu^{(2)} &=& {{\itbf W}^{(1)}}^T \boldsymbol{\Theta}^{(0)} {\itbf W}^{(1)} ,
\end{eqnarray}
and $ {\itbf W}^{(1)}$ is solution of $ \boldsymbol{\Theta}^{(0)}  {\itbf W}^{(1)} =( \boldsymbol{\Psi}^{(1)} -\mu^{(1)})
{\itbf W}^{(0)}$ and is orthogonal to $ {\itbf W}^{(0)}$.
If, for instance, $\Gamma_{jl}\equiv \Gamma >0$ for all $j\neq l$, then
$$
W_{jl}^{(1)}= - \frac{c_N}{\Gamma N} \big( \Lambda_j^{(1)}+ \Lambda_l^{(1)} -2\lambda^{(1)}\big) , \quad j\leq l,
$$
and
$$
\mu^{(2)}  =
\sum_{j\leq l} W_{jl}^{(1)} (  \boldsymbol{\Theta}^{(0)}  {\itbf W}^{(1)} )_{jl}= - \frac{2(N+2)}{N^2(N+1) \Gamma}  \sum_{j=1}^N (\Lambda_j^{(1)}-\lambda^{(1)})^2  .
$$
Note that
\begin{eqnarray}
\nonumber
\mu - 2\lambda  &=&  
\delta^2\big(\mu^{(2)}  -2 \lambda^{(2)}\big) +O(\delta^3)
\\
&=&
-\frac{2 \delta^2}{N^2(N+1) \Gamma }   \sum_{j=1}^N (\Lambda_j^{(1)}-\lambda^{(1)})^2 
+O(\delta^3) .
\label{eq:diffvp}
\end{eqnarray}
is negative-valued.

{\bf Exponential growth of the intensity fluctuations.}
It is a general feature that, for any matrix $\boldsymbol{\Gamma}$ and effective dissipation coefficients
$\Lambda_j$, we have
$\mu - 2\lambda  \leq 0$ (we have equality when there is no effective dissipation and this is a consequence of the forthcoming result (\ref{eq:si1}) and 
Cauchy-Schwarz inequality when there is dissipation).
The first two moments of the pointwise intensity 
$|\hat{p}(x,z)|^2$ for large $x$ are
\begin{align}
\EE [ |\hat{p}(x,z)|^2] =& \sum_{j=1}^N\frac{\phi_j(z)^2}{\beta_j} c_V V_j e^{- \lambda x} ,\\
\EE [ |\hat{p}(x,z)|^4 ] =& \sum_{j,l=1}^N \frac{\phi_j(z)^2\phi_l(z)^2}{\beta_j \beta_l} c_W W_{jl} e^{- \mu x} .
\end{align}
Without dissipation we have the following result for the fluctuations of the pointwise intensity:
$$ 
\frac{\EE[|\hat{p}(x,z)|^4]}{\EE[|\hat{p}(x,z)|^2]^2} \stackrel{x \to \infty}{\longrightarrow} 
\frac{2N}{N+1}  ,
$$
which is equal to $2$ when $N \gg 1$,
and with dissipation
\begin{equation}
\frac{\EE[|\hat{p}(x,z)|^4]}{\EE[|\hat{p}(x,z)|^2]^2} \stackrel{x \to \infty}{\simeq} 
\frac{2N}{N+1} 
\exp\big( - (\mu - 2\lambda) x \big)  ,
\label{eq:si1}
\end{equation}
that grows exponentially with the propagation distance (for very long distances, however, as $\mu - 2\lambda $
is very small as shown above).
Eq.~(\ref{eq:diffvp}) gives the expression of the exponential growth rate 
when dissipation is weak. In this regime the growth rate increases when the effective modal dissipation coefficients
become different from each other and decreases when the number of modes increases.

\section{Inverse problem}
\label{sec:ip}%
The goal of this section is to show that it is possible to estimate the statistical properties of the index of refraction
in the water column and the sea bottom properties from the sound pressure recorded by a vertical hydrophone array
and transmitted by distant time-harmonic sources.
This inverse problem can be formulated as a minimization problem that tries to match empirical quantities 
with theoretical ones that depend on the parameters to be estimated.
The theoretical model of the previous section shows that the correlation function of the sound pressure
has a non-trivial behavior that makes it possible to identify many relevant parameters as we explain below.
We first present the ALMA experimental data, then formulate the inverse problem, and finally estimate the model parameters.

\subsection{ALMA 2016 experiment}
\begin{figure}[!h]
  \centering
  \includegraphics[scale=0.25]{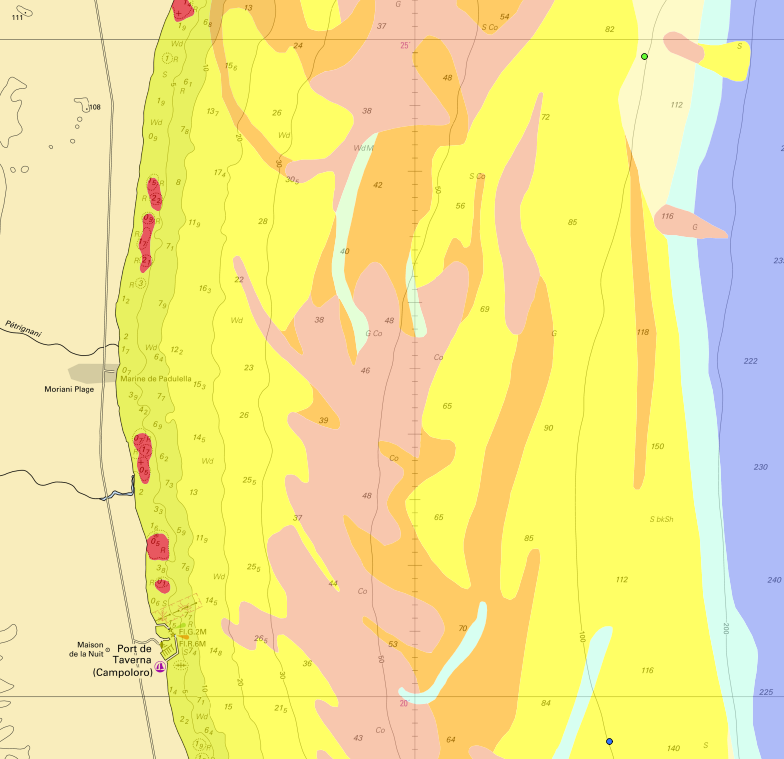}
  \caption{Top view of the experimental setup measurements region (source: data.shom.fr); blue dot (bottom right): source ($42^o 19.656$N $9^o 37.004$E); green dot (top right): receiver array ($42^o 24.693$N $9^o 37.363$E); yellow: sand sediment; orange: gravelly-sand sediment.}
  \label{carteALMA2016}
\end{figure}

ALMA is a series of at-sea experiments carried out by French DGA Naval Systems \cite{alma15}. In 2016, 
the experiment took place in november 7th-17th near the shores of North East Corsica \cite{alma17c} (cf. Figure~\ref{carteALMA2016}). 
There was a pinger immerged at fifty meters, and nine kilometers away, a passive array with $128$~hydrophones ($4$~verticals arms of $32$~hydrophones, the arms are $0.5~m$ apart from each other, the hydrophones are $0.15~m$ apart from each other) immerged at sixty meters. The sea bottom is relatively flat on the line between the pinger and the array (between $100~m$ and $115~m$), and composed of sand and gravelly-sand. Sea was calm, with roughness height around $0.1~m$ (sea state 1).
The pinger transmitted, for about 2 hours and a half, a sequence consisting of several time-harmonic waves, with a repetition rate of three minutes.
The sequence was a train made of two-second time-harmonic waves at the following $K=6$ frequencies: $2 kHz$, $5~kHz$, $7~kHz$, $9~kHz$, $11~kHz$, $13~kHz$. Acquisition was sampled at $48~kHz$. The Fourier coefficients for each frequency were extracted using a Hann window function with a duration of $1~s$.

Each arm has the form of a vertical array at a distance $x_{\rm a}$ of the source immerged at $z_0$. The hydrophones are regularly spaced between the depths $z_m$ and $z_M$, with $0 < z_m < z_M < z_{\rm b} $.
At frequency $\omega$, for $y \in \left[ 0 , z_M - z_m \right]$, we denote the spatial correlation between hydrophones at a distance $y$ from each other by $C_{x_{\rm a}}(y)$ (cf. Eq.~(\ref{cst}) for the theoretical expression predicted by our model). 
Then we define the correlation radius $r$ as the 
half-width at half-maximum (i.e. the distance between hydrophones for which   the correlation is $1/2$):
$$
C_{x_{\rm a}}(r) = \frac{1}{2} C_{x_{\rm a}}(0)  =\frac{1}{2} .
$$
From the experimental data at frequencies $f_1, \dots , f_K$ (with $\omega=2\pi f$) we can extract the experimental correlation radii $r_{\rm e}(f_1), \dots , r_{\rm e}(f_K) $. For $\Phi$ a set of parameters of our theoretical model (assuming that we know the sound speed $c_{\rm w}$ and density $\rho_{\rm w}$ in water and the depth of the sea bottom),
$$
  \Phi = \left( c_{\rm s} , \rho_{\rm s} , \nu_{\rm s} , \sigma , \ell_{\rm v} , \ell_{\rm h} \right)   ,
$$
we can compute the theoretical correlation radii $r_{\rm t}(f_1,\Phi) , \dots , r_{\rm t}(f_K,\Phi)$. Then we can define the misfit function
\begin{equation}
\label{eq:optcost}
  E(\Phi) = \sum_{k=1}^K \left( r_{\rm t}(f_k , \Phi) - r_{\rm e}(f_k) \right)^2  ,
\end{equation}
and we can determine the parameters $\Phi$ by minimizing over $\Phi$ the misfit function $E(\Phi)$.

\subsection{Estimation of the model parameters}
\begin{figure}[h!]
  \centering
  \includegraphics[scale=0.25]{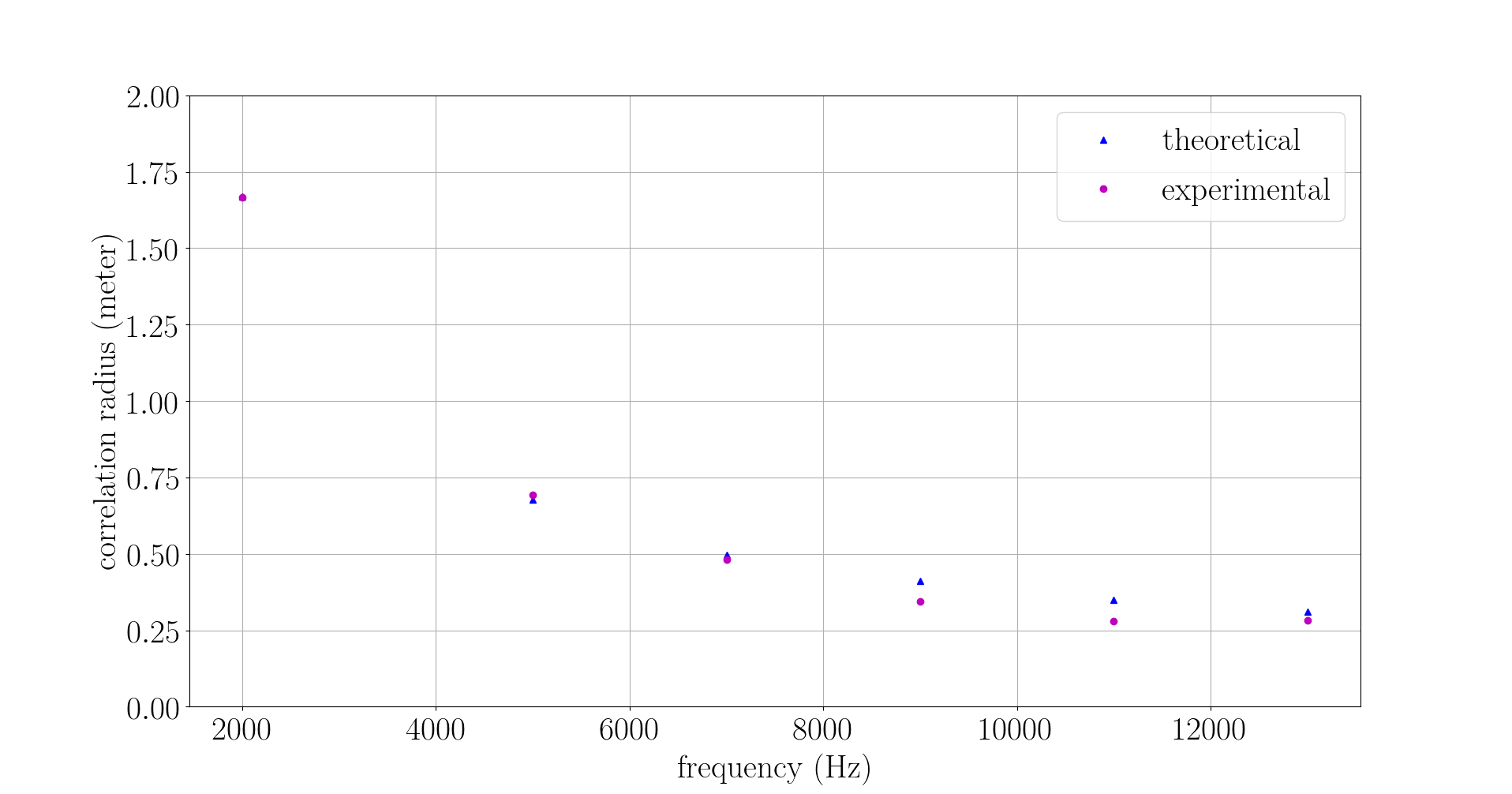}
  \caption{Comparisons between experimental and theoretical correlation radii. }
  \label{rc2016}
\end{figure}

\begin{figure}[h!]
\centering
\begin{tabular}{cc}
\includegraphics[scale=0.18]{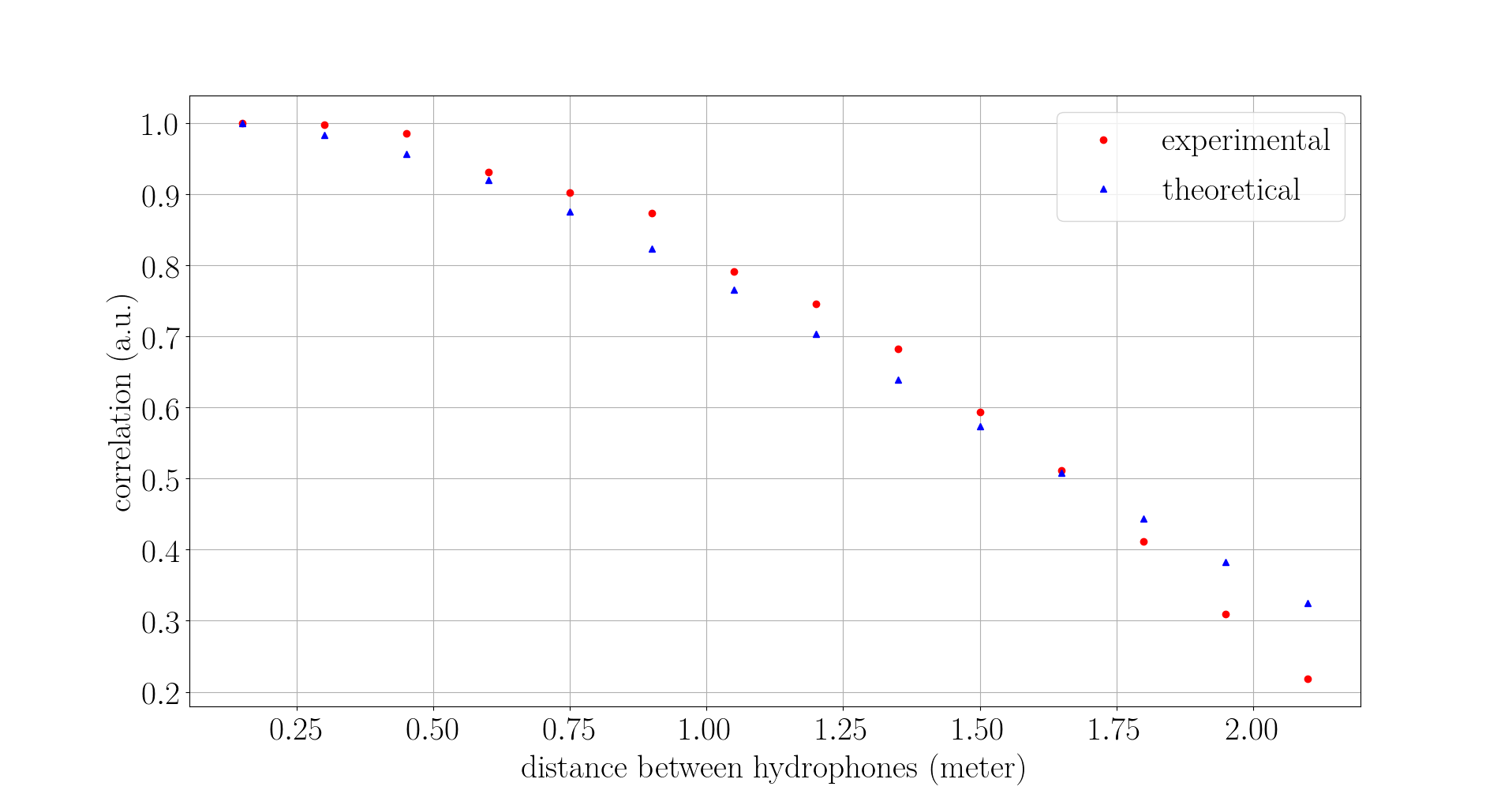} & \includegraphics[scale=0.18]{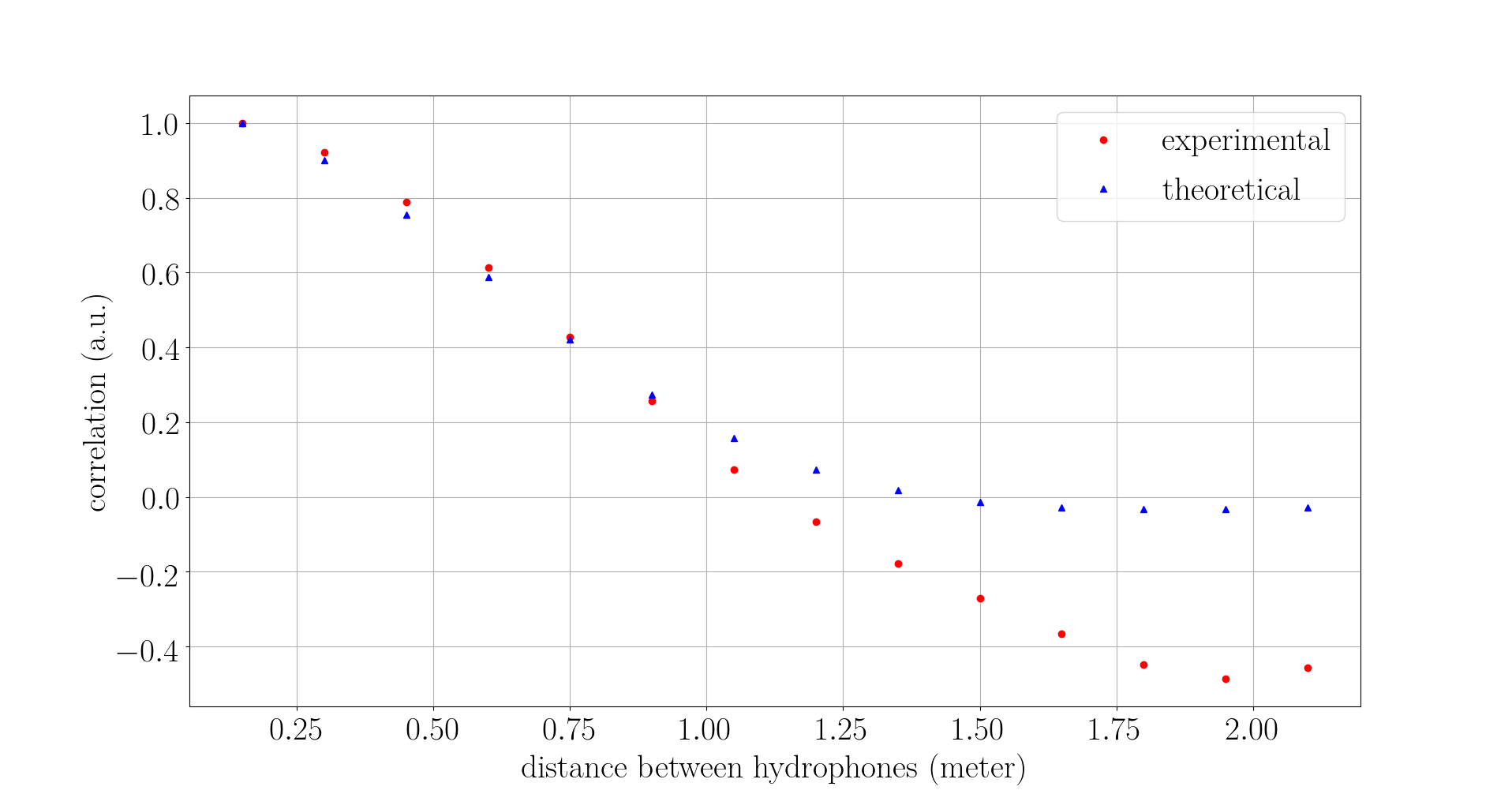}\\
$2~kHz$ & $5~kHz$\\
\includegraphics[scale=0.18]{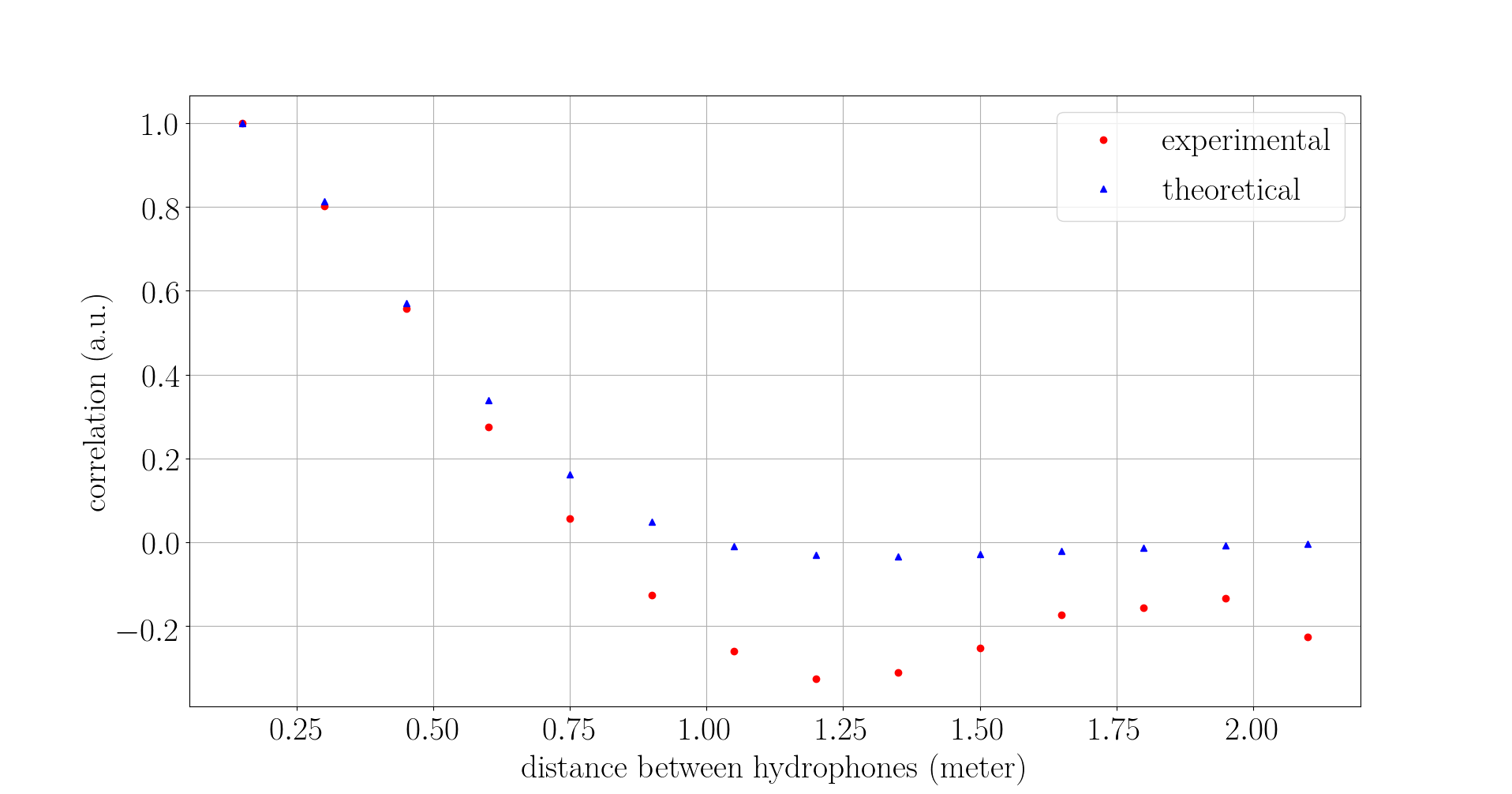} &  \includegraphics[scale=0.18]{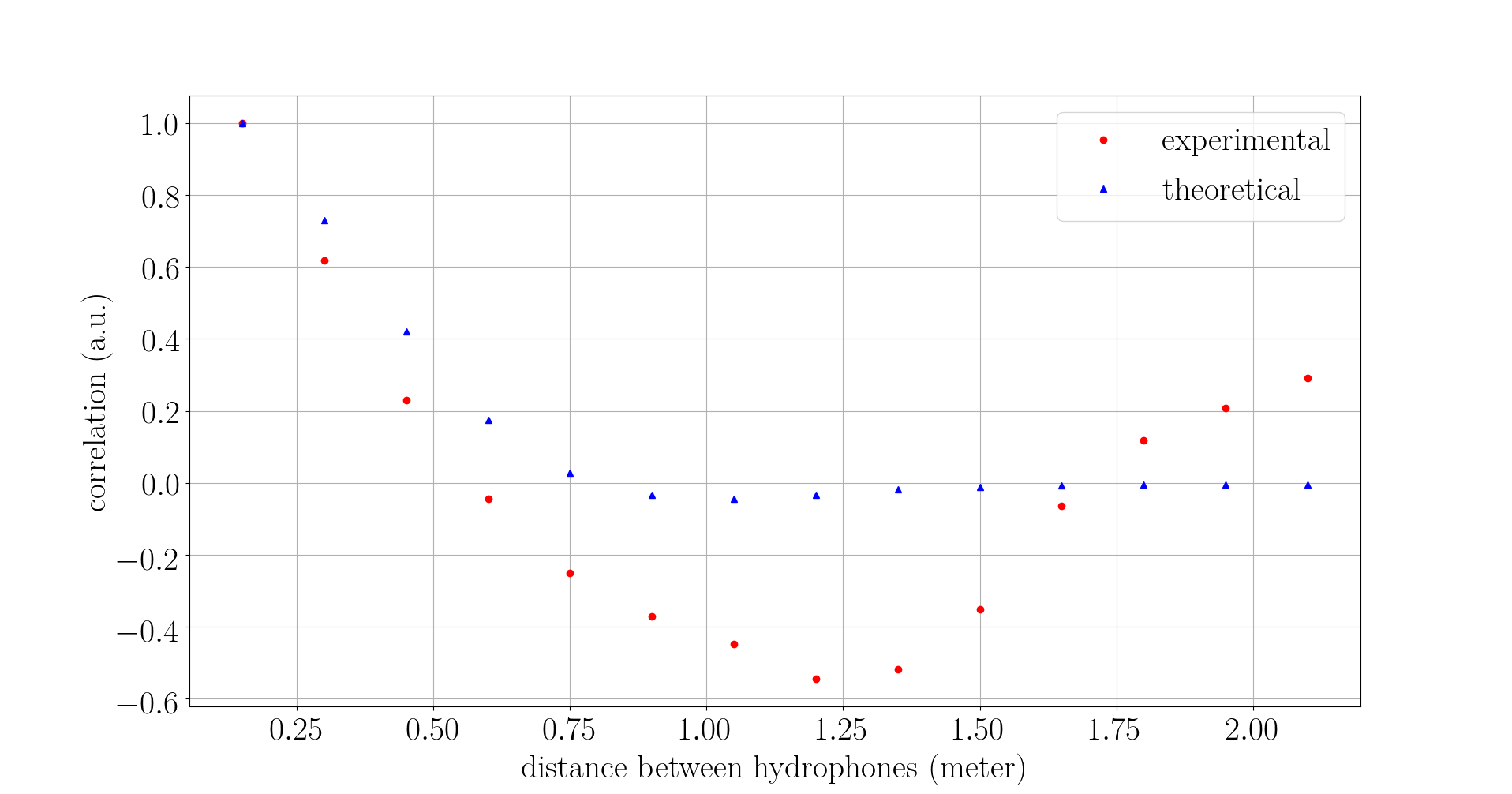}\\
$7~kHz$ & $9~kHz$\\
\includegraphics[scale=0.18]{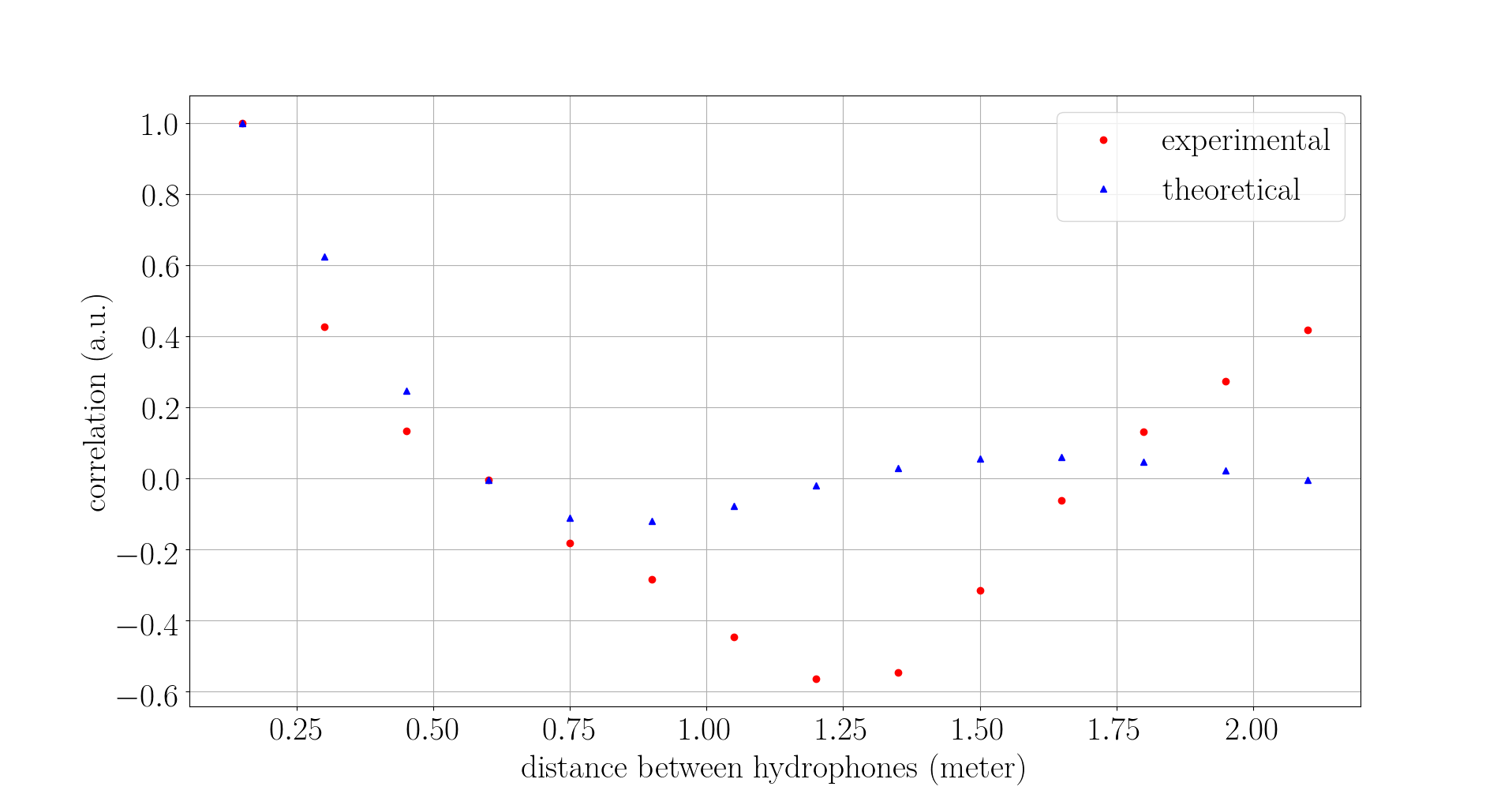} & \includegraphics[scale=0.18]{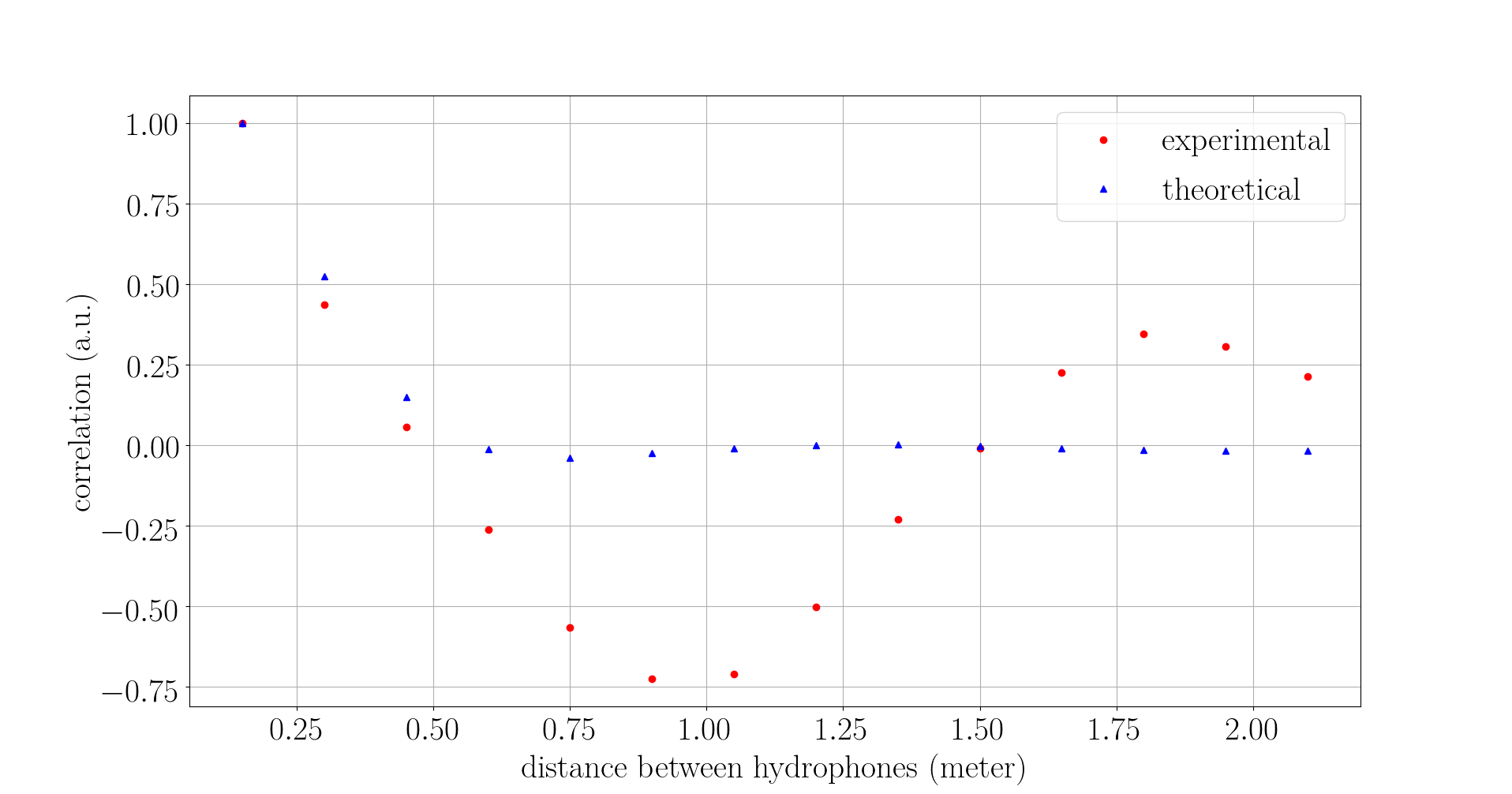}\\
$11~kHz$ & $13~kHz$
\end{tabular}
\caption{Experimental and theoretical correlation functions at different frequencies.}
\label{cf13Khz}
\end{figure}

We assume that sound speed is constant in water, taking value $c_{\rm w}=1523 m/s$ from CTD (conductivity, temperature, and depth) measurements, and we set density of water $\rho_{\rm w}$ at $1000~kg/m^3$ (we could have taken a more precise value, but density has little influence on the model).  Here, we have $x_{\rm a} = 9000~m$. Both source and array are in the middle of the water column.
Figure~\ref{rc2016} presents the correlation radii obtained from experimental data for each frequency, and the theoretical correlation radii 
predicted by our model with the parameters determined by minimization of the misfit function  (\ref{eq:optcost}):
    $ c_{\rm s} = 1630~m/s $;
    $\alpha = 1.09~dB / \textrm{wavelength} $;
    $\sigma = 0.002 $ (relative fluctuation);
    $\ell_{\rm v} = 30~m$;
    $\ell_{\rm h} = 100~m$;
    $\rho_{\rm s} = 1700~kg/m^3 $.
Parameter $\nu_{\rm s}$ is computed from the attenuation coefficient $\alpha$ expressed in decibel by wavelength. 
The values of the parameters seem compatible with experimental measurements carried out in similar environments \cite{rb04}.
Figure~\ref{cf13Khz}  shows
 theoretical and experimental correlation functions at the different frequencies.

\begin{figure}[h!]
  \centering
  \begin{tabular}{cc}
  \includegraphics[scale=0.18]{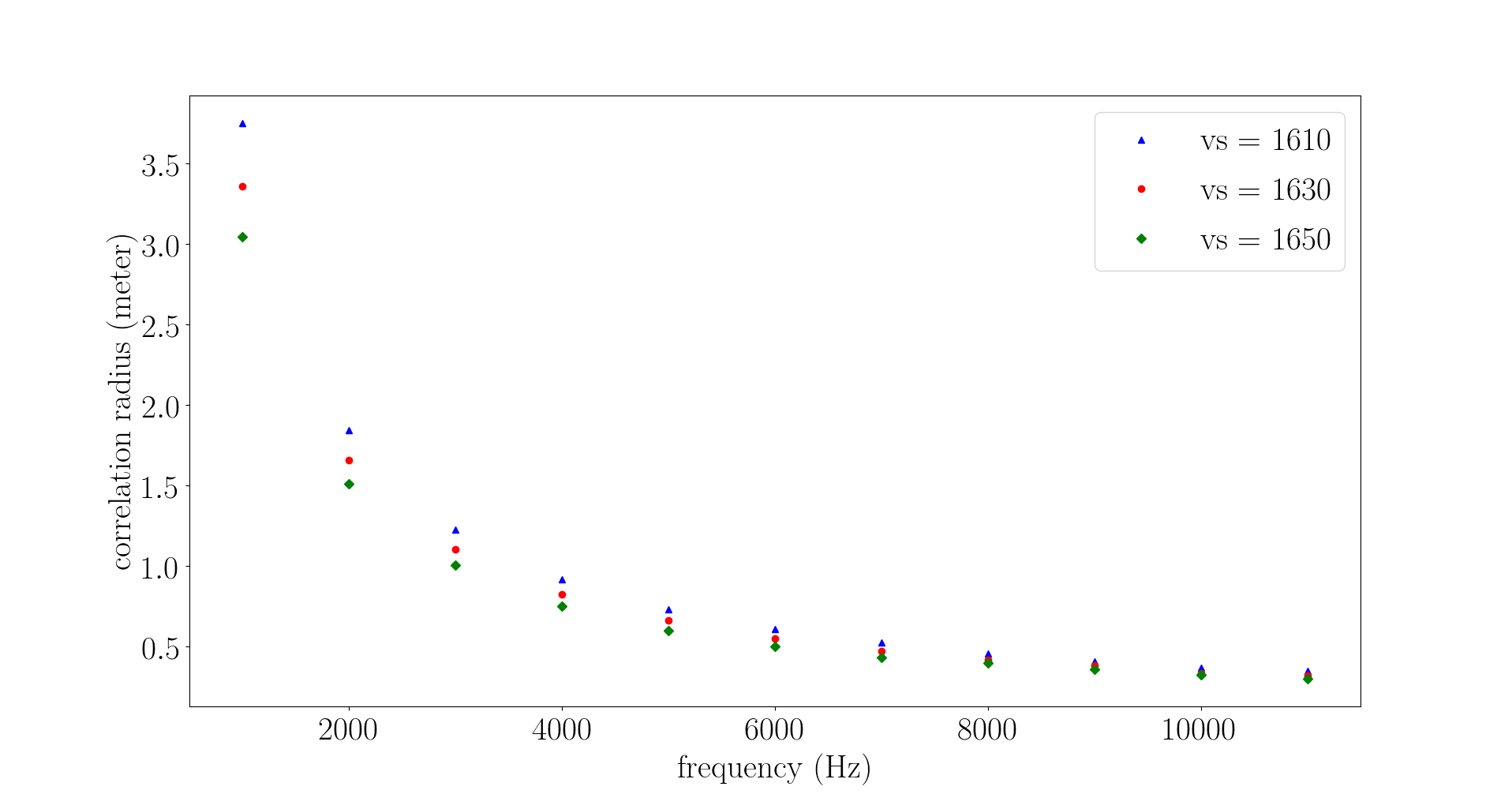} & \includegraphics[scale=0.18]{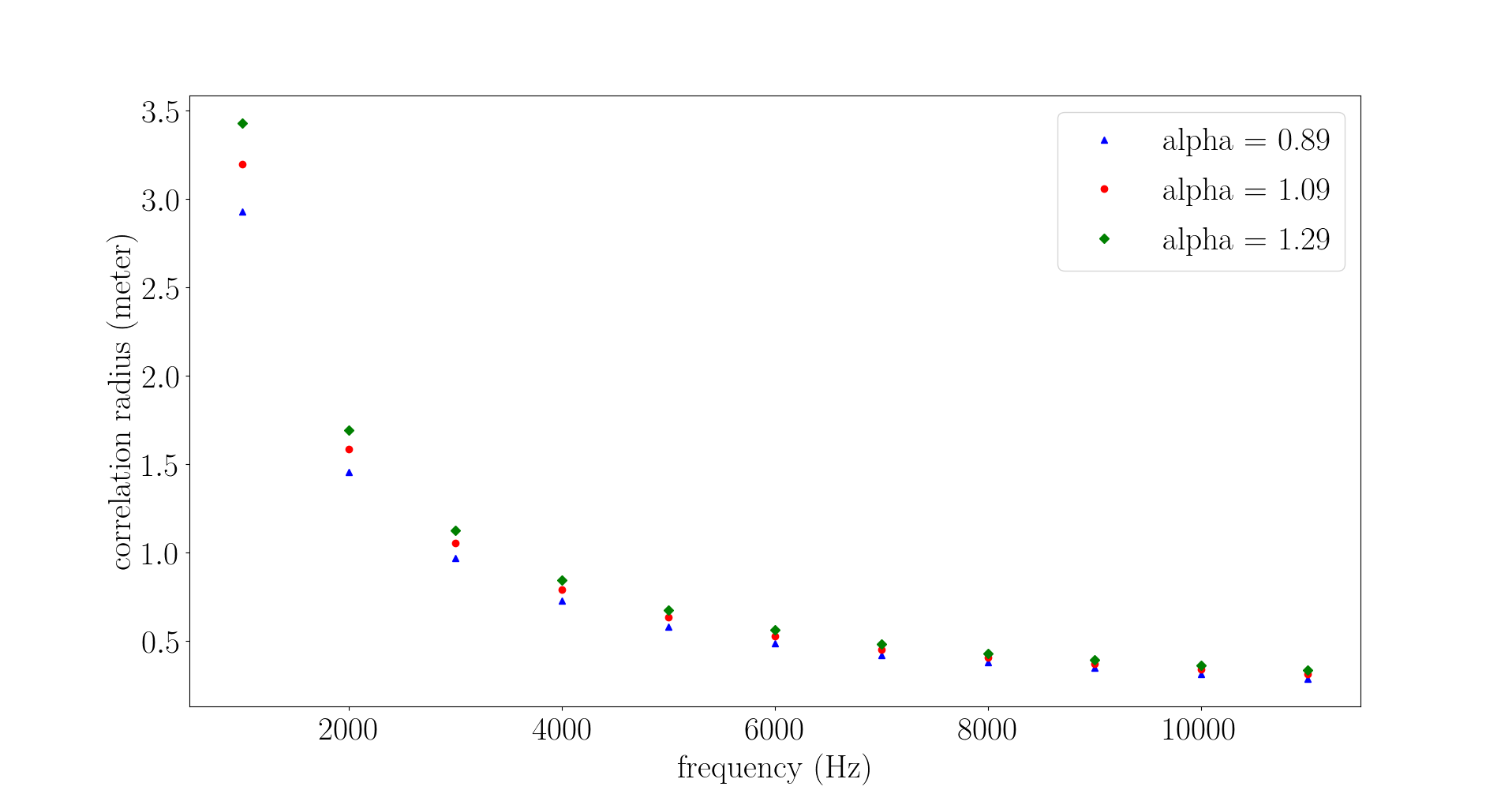}\\
    \includegraphics[scale=0.18]{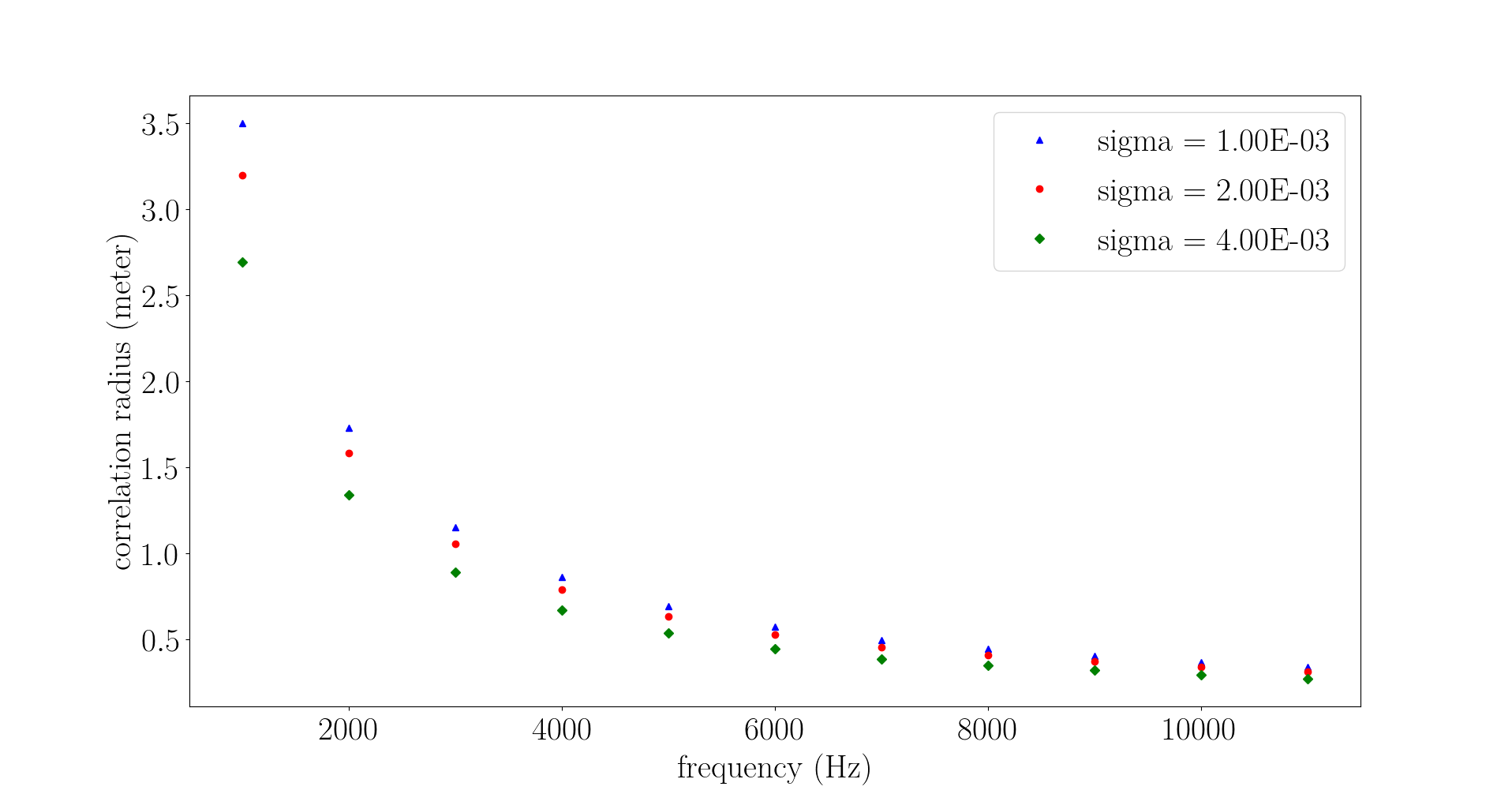}  & \includegraphics[scale=0.18]{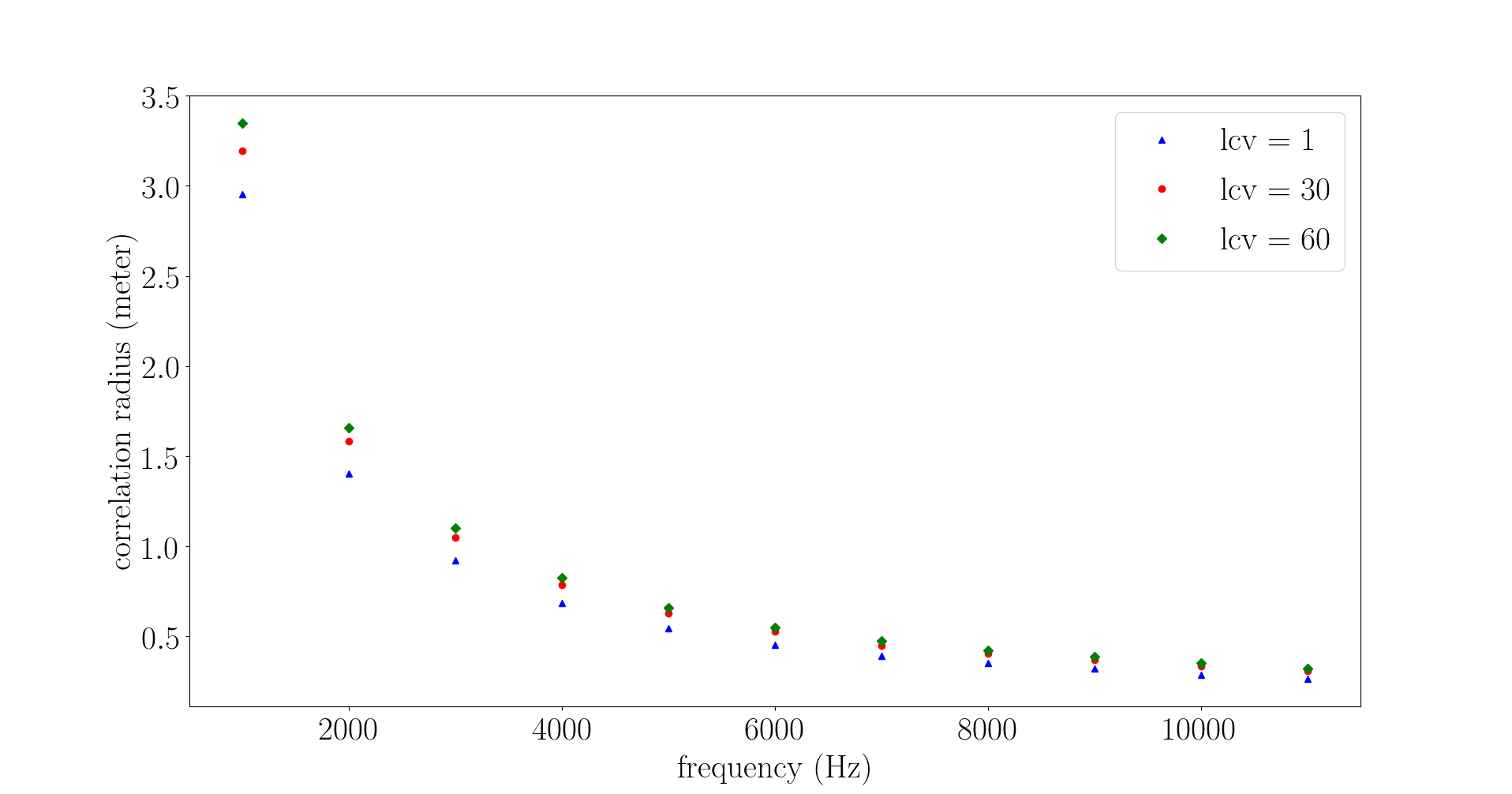}
    \end{tabular}
  \caption{Sensibility of the theoretical correlation radii to $c_{\rm s}$ (vs, in $m/s$, top left), $\alpha$ (alpha, in $dB / \textrm{wavelength} $, top right),
   $\sigma$ (sigma, bottom left), and $\ell_{\rm v}$ (lcv, in $m$, bottom right).}
  \label{AS_lcv}
\end{figure}

We carry out a basic (one-at-a-time) sensitivity analysis: we modify one by one the optimal values we have found, letting the others unchanged.
Results are in Figure~\ref{AS_lcv}. We see that $c_{\rm s}$ has a noticable influence, as a $20~m/s$ change in sediment sound speed significantly modifies the  correlation radii. $\sigma$ is also important, as dividing or multiplying it by two has a strong influence. $\alpha$ and $\ell_{\rm v}$ have less influence, but are still significant. Finally, for reasonable physical changes of $\rho_{\rm s}$ and $\ell_{\rm h}$, we do not see noticeable effects in the results, whatever the frequency is.

The inspection of the numerical values of the two terms that determine $\Lambda_j$ (cf. Eq.~(\ref{def:Lambdaj})) give interesting information.
We can see that for each frequency, the radiative terms are approximately ten times smaller than the attenuation ones. 
So we can claim that, in the experimental configuration addressed in this paper, radiation effects can be neglected.

Finally, if we denote by  $(z_n)_{1\leqslant n \leqslant N_{\rm h}}$ the depths of the $N_{\rm h}$ hydrophones of the array, 
we can compare the theoretical frequency-dependent scintillation indices defined by
\[
S_{\rm t} = \frac{1}{N_{\rm h}} \sum_{n=1}^{N_{\rm h}} \frac{\EE[|\hat{p}(x_{\rm a},z_n)|^4] - \EE[|\hat{p} (x_{\rm a},z_n)|^2]^2}{\EE[|\hat{p}(x_{\rm a},z_n)|]^2} ,
\]
with the experimental values $S_{\rm e} $ determined by empirical averages instead of expectations.
Table \ref{tab:ind} shows experimental scintillation indices compared to theoretical ones, with the optimal parameters $\Phi$ 
determined from the correlation radii. In the experimental and theoretical configurations scintillation indices are approximately equal to one. One should 
have longer propagation distances to observe scintillation indices larger than one, as already suggested in Ref.~~\onlinecite{creamer}.
\begin{table}
\[
  \begin{array}{|c|c|c|c|c|c|c|}
    \hline
  ~ &   \textrm{$2~kHz$} &   \textrm{$5~kHz$} &   \textrm{$7~kHz$} &   \textrm{$9~kHz$}  &   \textrm{$11~kHz$}  &   \textrm{$13~kHz$}\\     \hline
\textrm{experimental}  & 0.82 & 0.89 & 1.06 & 1.41 & 1.16 & 1.89 \\     \hline
  \textrm{theoretical} & 0.96 & 0.98 & 0.99& 0.99& 0.99& 0.99 \\ \hline
\end{array}
\]
\caption{Experimental and theoretical scintillation indices at different frequencies.}
\label{tab:ind}
\end{table}

\section{Conclusion}
This paper proposes a complete description of the statistics of the mode amplitudes of the sound pressure in a shallow-water waveguide.
The effective parameters (frequency- and mode-dependent attenuation, dispersion, and coupling) are identified from first principles
and expressed in terms of the statistical properties of the index of refraction of the water column and the sea bottom properties.
This theoretical analysis makes it possible to formulate an inverse problem for the estimation of the model parameters 
from the sound pressure recorded by a vertical  hydrophone array and transmitted by distant time-harmonic sources.
This inverse problem is solved using data collected during the ALMA 2016 experiment.

In the experimental configuration addressed in this paper, it turns out that radiation effects can be neglected (compared to dissipation in the sediments) 
and non-Gaussian scintillation effects can be neglected as well.
The extraction of the  frequency-dependent correlation radii of the recorded sound pressure signals 
makes it possible to estimate different acoustic and geoacoustic parameters:
The sediment sound speed $c_{\rm s}$ and the standard deviation of the index of refraction in the water column $\sigma$ can be robustly estimated 
because they have strong effects on the correlation radii of the sound pressure. The attenuation parameter in the sediments $\alpha$ and the vertical  correlation radius of the index of refraction in the water column $\ell_{\rm v}$ can be estimated, because they have less influence, but are still significant. 
The sediment density $\rho_{\rm s}$ and the horizontal correlation radius $\ell_{\rm h}$ are difficult to estimate because they 
have small effects of the correlation radii of the sound pressure.

\section*{Acknowledgements}
We thank Dominique Fattaccioli and Gaultier Real from DGA Naval Systems
for interesting and stimulating discussions and for making us available the ALMA 2016 data.

\appendix
\section{The effective Markovian system for the complex mode amplitudes}
\label{app:A}
We assume that the reduced wavenumbers $\beta_{j}$  are distinct, for $1 \le j \le N$. Then, for any $x_\infty>0$, the  process 
$
\big( (\hat{a}_{j}^\eps(x) )_{j=1}^{N}, (\hat{a}_{\gamma}^\eps(x) )_{\gamma \in (0,k_{\rm s}^2)} \big) 
$
converges in distribution in ${\cal C}^0([0,x_\infty], \mathbb{C}^{N}  \times L^2(0,k_{\rm s}^2) )$, where $\mathbb{C}^{N}  \times L^2(0,k_{\rm s}^2)$
is equipped with the weak topology,  to the Markov process  
$
\big( ({a}_{j}(x) )_{j=1}^{N} , ({a}_{\gamma} (x))_{\gamma \in (0,k_{\rm s}^2)} \big)
$
with infinitesimal generator
$
{\cal L}= {\cal L}^1+{\cal L}^2+{\cal L}^3 ,
$
where ${\cal L}^j$, $1 \leq j \leq 3$, are the differential operators:
\begin{align}
\nonumber {\cal L}^1 = & \frac{1}{2} 
\sum_{j,l=1}^{N} 
\Gamma_{jl} \big( 
{a}_{j} \overline{{a}_{j}} \partial_{{a}_{l}} \partial_{\overline{{a}_{l}}}
+
{a}_{l} \overline{{a}_{l}} \partial_{{a}_{j}} \partial_{\overline{{a}_{j}}}
-
{a}_{j}  {a}_{l} \partial_{{a}_{j}} \partial_{ {a}_{l}}
-
\overline{{a}_{j}}  \overline{{a}_{l}} \partial_{\overline{{a}_{j}}} \partial_{ \overline{{a}_{l}}}\big)
{\bf 1}_{ j\neq l}
\\
\nonumber &+
 \frac{1}{2} 
\sum_{j,l=1}^{N} 
\Gamma^1_{j l} \big( 
{a}_{j} \overline{{a}_{l}} \partial_{{a}_{j}} \partial_{\overline{{a}_{l}}}
+
\overline{{a}_{j}} {a}_{l} \partial_{\overline{{a}_{j}}} \partial_{{a}_{l}}
{a}_{j}  {a}_{l} \partial_{{a}_{j}} \partial_{ {a}_{l}}
-
\overline{{a}_{j}}  \overline{{a}_{l}} \partial_{\overline{{a}_{j}}} \partial_{ \overline{{a}_{l}}}\big)
\\
&+ \frac{1}{2} 
\sum_{j=1}^{N} \big( \Gamma_{jj} - \Gamma^1_{jj}\big)
\big( {a}_{j} \partial_{{a}_{j}} + \overline{{a}_{j}} \partial_{\overline{{a}_{j}}}
\big)
+\frac{i}{2} 
\sum_{j=1}^{N}  \Gamma^{s}_{jj}  
\big( {a}_{j} \partial_{{a}_{j}} - \overline{{a}_{j}} \partial_{\overline{{a}_{j}}}
\big)  , \label{eq:defL1}
\\
{\cal L}^2
=&
-\frac{1}{2}  
\sum_{j=1}^{N} ( \Lambda_{j}  +i\Lambda^{s}_{j} ) {a}_{j} \partial_{{a}_{j}} 
+
( \Lambda_{j}  - i\Lambda^{s}_{j} )  \overline{{a}_{j}} \partial_{\overline{{a}_{j}}}  ,
 \label{eq:defL2}\\{\cal L}^3
=&
i  
\sum_{j=1}^{N}  \kappa_{j}   \big( {a}_{j} \partial_{{a}_{j}} 
- \overline{{a}_{j}} \partial_{\overline{{a}_{j}}} \big) .
 \label{eq:defL3}
\end{align}
In these definitions we use the  classical complex derivative:
if $ \zeta=\zeta_r+i\zeta_i$, then $\partial_\zeta=(1/2)(\partial_{\zeta_r}-i \partial_{\zeta_i})$ and
$\partial_{\overline{\zeta}} =(1/2)(\partial_{\zeta_r} +i \partial_{\zeta_i})$,
and the coefficients of the operators (\ref{eq:defL1}-\ref{eq:defL3}) are defined for 
indices $j ,l= 1, \ldots, N$, as follows: 

 \noindent For all $j \neq l$, $\Gamma_{jl}$ is defined by (\ref{def:Gammalj}) and 
\begin{align*}
\Gamma^{s}_{j l} =&
\frac{\omega^4}{2 \beta_{j}\beta_{l}}
\int_0^\infty \EE\big[ C_{jl}^{\rm w}(0) C_{jl}^{\rm w}(x) \big] \sin \big[(\beta_{l}-\beta_{j})x\big] dx .
\end{align*}
For all $j,l$:
\begin{align*}
\Gamma^1_{jl} =&
\frac{\omega^4}{2 \beta_{j}\beta_{l}}
\int_0^\infty \EE\big[ C_{jj}^{\rm w}(0) C_{ll}^{\rm w}(x) \big]   dx.
\end{align*}
For all $j$, $\Lambda_j$ is defined by (\ref{def:Lambdaj}) and
\begin{align*}
\Gamma_{jj} =& -\hspace{-0.05in}\sum_{l =1,l\neq j}^{N} \Gamma_{jl} ,\quad \quad
\Gamma^{s}_{jj} = -\hspace{-0.05in}\sum_{l =1,l\neq j}^{N} \Gamma^{s}_{jl}  ,\\
\Lambda_{j}^{s} =&\hspace{-0.05in}    \int_0^{k_{\rm s}^2} \hspace{-0.05in}\frac{\omega^4}{2\sqrt{\gamma}\beta_{j}}
\int_0^\infty \EE \big[ C_{j \gamma}^{\rm w}(0) C_{j\gamma}^{\rm w}(x) \big] 
\sin \big[ (\sqrt{\gamma}-\beta_{j})x\big] dx d\gamma, \\
\kappa_{j} = &  \hspace{-0.05in}   \int_{-\infty}^0  \frac{\omega^4}{2\sqrt{|\gamma|}\beta_{j}}
\int_0^\infty \EE \big[ C_{j \gamma}^{\rm w}(0) C_{j \gamma}^{\rm w}(x) \big] \cos ( \beta_{j} x) e^{-\sqrt{|\gamma|}x}dx d\gamma.
\end{align*} 

Remarks: \\
1) The convergence result holds in the weak topology. This means that we can only compute quantities of the 
form $\EE [ F( a_1,\ldots,a_N, \int_0^{k_{\rm s}^2} \alpha_\gamma a_\gamma d\gamma)]$
for any test function $\alpha \in L^2(0,k_{\rm s}^2)$ and $F:\mathbb{R}^{N+1}\to \mathbb{R}$,
which are the limits of $\EE [ F( \hat{a}^\eps_1,\ldots,\hat{a}^\eps_N, \int_0^{k_{\rm s}^2} \alpha_\gamma \hat{a}^\eps_\gamma d\gamma)]$ as $\eps \to 0$. \\
2) The generator ${\cal L}$ does not involve $\partial_{a_{\gamma}}$. Therefore 
$  (\hat{a}_{j}^\eps (x))_{j=1}^{N}   $
converges in distribution in ${\cal C}^0([0,x_\infty], \mathbb{C}^{N}  )$ to
the Markov process  $  ({a}_{j}(x) )_{j=1}^{N}   $
with  generator ${\cal L}$. The weak and strong topologies are the same  in $\mathbb{C}^N$,
so we can compute any moment of the 
form $\EE [ F( a_1,\ldots,a_N)]$, which are the limits of $\EE [ F( \hat{a}^\eps_1,\ldots,\hat{a}^\eps_N)]$.\\
3) ${\cal L}_1$ is the contribution of the coupling between guided modes, 
which gives rise to power exchange between the guided modes;
${\cal L}_2$  is the contribution of the coupling between guided and radiating modes, 
which gives rise to power leakage from the guided modes to the radiating ones (effective diffusion)
and addition of frequency-dependent phases on the guided modes (effective dispersion);
${\cal L}_2$ also contains the effective mode-dependent term due to dissipation in the sediments;
${\cal L}_3$  is the contribution of the coupling between guided and evanescent modes,
which gives rise to additional phase terms on the guided modes (effective dispersion \cite{garnier_evan}).\\
4) If the generator ${\cal L}$ is applied to a test function that depends only on the mode powers $  (P_{j} = |{a}_{j}|^2 )_{j=1}^{N}   $,
then the result is a function that depends only on $ (P_{j} )_{j=1}^{N}  $. Thus, 
the mode powers $  (P_{j}(x) )_{j=1}^{N}   $ define  a Markov process, with infinitesimal generator defined by (\ref{eq:genLP}).\\
5) The radiation mode amplitudes  remain constant on $L^2(0,k_{\rm s}^2)$, equipped with the weak topology, as $\eps \to 0$. However, this does not
describe the power $  \int_0^{k_{\rm s}^2} |\hat{a}^\eps_{\gamma}|^2 d\gamma$ transported by the radiation modes,  because the convergence does not hold in the strong topology of $L^2(0,k_{\rm s}^2)$ so we do not have $ \int_0^{k_{\rm s}^2} |\hat{a}^\eps_{\gamma}|^2 d\gamma \to  \int_0^{k_{\rm s}^2} |a_{\gamma}|^2 d\gamma$ as $\eps \to 0$.


\begin{thebibliography}{99}

\bibitem{beran}
M. J. Beran and S. Frankenthal,
Volume scattering in a shallow channel,
J. Acoust. Soc. Am. {\bf 91}, 3203-3211 (1992).

\bibitem{borcea15}
L. Borcea, J. Garnier, and C. Tsogka, A quantitative study of source imaging in random waveguides, 
Commun. Math. Sci. {\bf 13}, 749-776 (2015).

\bibitem{colosi09}
J. A. Colosi and A. Morozov,  
Statistics of normal mode amplitudes in an ocean with random sound speed perturbations: 
Cross mode coherence and mean intensity,
J. Acoust. Soc. Am. {\bf 126}, 1026-1035 (2009).

\bibitem{colosi12}
J. A. Colosi, T. F. Duda, and A. K. Morozov,
Statistics of low-frequency normal-mode amplitudes in an ocean with random sound-speed perturbations: 
Shallow-water environments,
J. Acoust. Soc. Am. {\bf 131}, 1749-1761 (2012).

\bibitem{creamer}
D. Creamer, 
Scintillating shallow water waveguides,
J. Acoust. Soc. Am. {\bf 99}, 2825-2838 (1996).

\bibitem{dozier}
L. B.  Dozier and F. D.  Tappert, 
Statistics of normal-mode am- plitudes in a random ocean. I. Theory, 
J. Acoust. Soc. Am. {\bf 63}, 353-365 (1978). 

\bibitem{dozier2}
L. B.  Dozier and F. D.  Tappert, 
Statistics of normal-mode am- plitudes in a random ocean. II. Computations,
J. Acoust. Soc. Am. {\bf 64}, 533-547 (1978).

\bibitem{dozier83}
L. B. Dozier, 
A coupled mode model for spatial coherence of bottom-interacting energy,
in Proceedings of the Stochastic Modeling Workshop, 
edited by C. W. Spofford and J. M. Haynes, ARL-University of Texas, Austin TX, 1983.

\bibitem{alma15}
D. Fattaccioli,
ALMA: a new experimental acoustic system to explore coastal and shallow waters, 
UACE 2015 Proceedings (Chania, Greece).

\bibitem{alma17c}
D. Fattaccioli and G. Real,
The DGA "ALMA" Project: an overview of the recent improvements of the system capabilities and of the at-sea campaign ALMA-2016,
UACE 2017 Proceedings (Skiathos, Greece),
available at www.uaconferences.org/docs/2017\_papers/649\_UACE2017.pdf

\bibitem{flatte}
S. M.  Flatt\'e, R. Dashen,  W. H. Munk, K. M. Watson, and F. Zachariasen, 
\emph{Sound Transmission Through a Fluctuating Ocean},
Cambridge University Press, Cambridge, 1979. 

\bibitem{FGPSbook}
J.-P. Fouque, J. Garnier, G. Papanicolaou, and K. S\o lna,
\emph{Wave Propagation and Time Reversal in Randomly Layered Media},
Springer, New York, 2007

\bibitem{garnier_evan}
J. Garnier, 
The role of evanescent modes in randomly perturbed single-mode waveguides,
Discrete and Continuous Dynamical Systems B {\bf 8}, 455--472  (2007). 

\bibitem{GP07}
J. Garnier and G. Papanicolaou, 
Pulse propagation and time reversal in random waveguides,
SIAM J. Appl.  Math. {\bf 67},  1718--1739 (2007).

\bibitem{garrett}
C. Garrett  and W. H. Munk, 
Internal waves in the ocean,
Annu. Rev. Fluid Mech. {\bf 11}, 339--369 (1979).

\bibitem{gomez}
C. Gomez,
Wave propagation in shallow-water acoustic random waveguides,
Commun. Math. Sci. {\bf 9},  81--125 (2011).

\bibitem{jensen}
F. B. Jensen, W. A. Kuperman, M. B. Porter, and H. Schmidt, 
{\it Computational Ocean Acoustics},
Springer, New York, 1993. Chap. 5.

\bibitem{kohler77}
W. Kohler and G. Papanicolaou,
Wave propagation in  randomly inhomogeneous ocean,
in Lecture Notes in Physics, Vol. 70,
J. B. Keller and J. S. Papadakis, eds.,
Wave Propagation and Underwater Acoustics,
Springer Verlag, Berlin, 1977.


\bibitem{rb04}
M. D. Richardson and K. B. Briggs,
Empirical predictions of seafloor properties based on remotely measured sediment impedance,
AIP Conference Proceedings {\bf 728}, 12--21 (2004).


\bibitem{wilcox}
C. Wilcox, Spectral analysis of the Pekeris operator in the theory of acoustic wave propagation
in shallow water, Arch. Rational Mech. Anal. {\bf 60}, 259--300 (1976).


\end{thebibliography}
\end{document}